\documentclass[11pt]{article}
\usepackage{amsmath,amssymb,amsfonts}
\usepackage{geometry}
\usepackage{hyperref}
\usepackage{algorithm}
\usepackage{graphicx}
\usepackage{algpseudocode}
\usepackage{authblk}

\newcommand{\tu}{\tilde{u}}
\newcommand{\cL}{\mathcal{L}}
\newcommand{\ueff}{u_\text{eff}^\ast}

\title{On-the-Fly Lifting of Coarse Reaction-Coordinate Paths to Full-Dimensional Transition Path Ensembles}

\author[1,2]{Christof Schuette}
\author[1]{Alexander Sikorski}
\author[1]{Jakob Kresse}
\author[1]{Marcus Weber}

\affil[1]{Zuse Institute Berlin, Berlin, Germany}
\affil[2]{Freie Universität Berlin, Institute of Mathematics and Computer Science, Berlin, Germany}

\date{\today}

\begin{document}
\maketitle

\begin{abstract}
Effective dynamics on a low-dimensional collective-variable (CV) or latent space can be simulated far more cheaply than the underlying high-dimensional stochastic system, but exploiting such coarse predictions requires lifting: turning a coarse CV trajectory into dynamically consistent full-dimensional states and path ensembles, without relying on global sampling of invariant or conditional fiber measures. We present a local, on-the-fly lifting strategy based on guided full-system trajectories. First an effective model in CV space is used to obtain a coarse reference trajectory. Then, an ensemble of full-dimensional trajectories is generated from a guided version of the original dynamics, where the guidance steers the trajectory to track the CV reference path. Because guidance biases the path distribution, we correct it via pathwise Girsanov reweighting, yielding a correct-by-construction importance-sampling approximation of the conditional law of the uncontrolled dynamics. We further connect the approach to stochastic optimal control, clarifying how coarse models can inform variance-reducing guidance for rare-event quantities. Numerical experiments demonstrate that inexpensive coarse transition paths can be converted into realistic full-system transition pathways (including barrier crossings and detours) and can accelerate estimation of transition pathways and statistics while providing minimal bias through weighted ensembles. 
\end{abstract}

\section{Introduction}
High-dimensional stochastic dynamics are ubiquitous in molecular simulation and related areas; a prototypical example for such dynamics is given by diffusion systems of the form
\begin{equation}
  dX_t \;=\; b(X_t)\,dt + \sigma\, dW_t, 
  \qquad X_t \in \mathcal{X}\subset\mathbb{R}^d,
  \label{eq:intro_full_sde}
\end{equation} 
or similar forms like Langevin dynamics.
In many applications, the dynamics is
\emph{metastable}: trajectories spend long times near a few dominant regions and only rarely
transition between them. Direct numerical time-stepping in full dimension is then expensive,
since small timesteps are required for stability and accuracy, while the rare transitions determine
the relevant long-time statistics and mechanisms. Transfer-operator and variational perspectives
provide a principled route to analyze such slow processes and to build reduced kinetic models
such as Markov state models (MSMs); see, e.g., \cite{SchuetteKlusHartmann2023,SchuetteSarich2014,BowmanPandeNoe2013,KlusNueskeKoltaiWuKevrekidisSchuetteNoe2018}.
One of the main challenges is to reliably find the transition regions and construct representative transition paths between the few dominant regions. 

A standard route to reduce complexity is to introduce a low-dimensional \emph{collective variable}
(CV) map $\xi:\mathcal{X}\to\mathbb{R}^m$, $m\ll d$, introducing a coarse representation of the full system, 
$Z_t := \xi(X_t)$, in a $m$-dimensional latent space. Among all the approaches to finding such CVs or \emph{reaction coordinates}, some specific \emph{dynamically informative} CVs are such that $Z_t$ captures the essential slow degrees of freedom \cite{ZhangSchuette2025}. Once an informative CV is
available, one may attempt to describe the evolution of $Z_t$ by an \emph{effective (coarse) dynamics}
in $\mathbb{R}^m$, which can often be simulated at substantially reduced cost. The construction,
analysis, and data-driven identification of such CVs and effective models has seen substantial recent
progress. Representative examples include time-lagged independent component analysis (TICA)
\cite{PerezHernandez2013}, variational approaches for Markov processes (VAMP) and their neural
implementations (VAMPnets) \cite{MardtPasqualiWuNoe2018}, reaction-coordinate flow ideas
\cite{WuNoe2024}, as well as Koopman/dominant-subspace learning approaches such as ISOKANN
\cite{RabbenRayWeber2020ISOKANN}. These algorithmic developments connect to theoretical
work on optimal reaction coordinates and exact dynamical coarse-graining without time-scale
separation \cite{LuVandenEijnden2014}, and to quantitative analyses of effective dynamics and
pathwise error bounds \cite{ZhangHartmannSchuette2016,LegollLelievreOlla2017,LelievreZhang2019}.
For related perspectives on choosing CVs by minimizing the deviation between effective and full
dynamics and on spectral properties of effective dynamics from conditional expectations, see
\cite{ZhangSchuette2025,NueskeKoltaiBoninsegnaClementi2021}.

Our central goal is to use knowledge about CVs and the related effective dynamics for computing transition rates and pathways of the full system under consideration in diffusive and molecular processes. There is extensive literature on computing transition kinetics and paths, including transition
path sampling (TPS) \cite{BolhuisChandlerDellagoGeissler2002TPS}, interface-based approaches such
as milestoning \cite{FaradjianElber2004Milestoning} and transition interface sampling (TIS)
\cite{MoroniVanErpBolhuis2004TIS}, and transition path theory (TPT)
\cite{EVE2006JSP,EVE2010TPTreview,MetznerSchuetteVandenEijnden2006TPTexamples}, for a conceptual review see
\cite{BerezhkovskiiSzabo2019CommittorsMilestones}. Some of these techniques like milestoning or TIS really make use of knowledge about CVs. However, all of these techniques suffer from the fact that they always explore the full-dimensional system and are not guided by a dynamically reliable low-dimensional
dynamics as given, e.g., by the effective dynamics. 

\paragraph{From coarse prediction to full-dimensional realization: the lifting bottleneck.}
Even when an effective model for $Z_t$ is available, a fundamental difficulty remains:
\emph{lifting}. Concretely, suppose the full dynamics \eqref{eq:intro_full_sde} has been simulated up to
time $T$, yielding a state $X_T$, and that the effective dynamics has been advanced further to a later
time $T_{\mathrm{new}}$, producing a coarse path
$\{Z_t\}$ over $[T,T_{\mathrm{new}}]$.
The lifting problem asks:

\begin{quote}
\emph{How can we generate a dynamically consistent full-dimensional path $\{X_t\}$ over $[T,T_{\mathrm{new}}]$ conditioned on the coarse path $\{Z_t\}$? In particular, if $\{Z_t\}$ is a transition path of the effective dynamics, how can we guarantee that $\{X_t\}$ is a representative transition path of the full-dimensional system? }
\end{quote}

\paragraph{Main idea: guided trajectories driven by a cheap coarse reference path.}
This work develops a \emph{local, on-the-fly} lifting strategy that uses the effective dynamics
to provide a coarse \emph{reference path} in CV space and then reconstructs full-dimensional
transition paths by simulating \emph{guided trajectories}. The key mechanism is to introduce
a controlled process
\begin{equation}
  dX_t^{u} \;=\; \bigl(b(X_t^{u}) +\sigma u(t,X_t^{u})\bigr)\,dt + \sigma\, dW_t,
  \label{eq:intro_controlled_sde}
\end{equation}
where the feedback control $u$ drives the full-system dynamics in a way that is driven by the CV $\xi$ and its effective dynamics on the latent space. Herein, we consider control forces of the general form
\begin{equation}
  u(t,x) \;=\; J_\xi(x)^\top \nabla \Psi(\xi(x),Z_t),
  \label{eq:general_control}
\end{equation}
with $J_\xi$ denoting the Jacobian of $\xi$ and  where $Z_t$ may be additional information on the latent space. As we will see, the choice of the function $\Psi$ depends on what problem is to be solved: (1) \emph{Fast exploration}, where one aims at finding previously unknown part of state space like relevant transition channels of the system, or (2) \emph{accurate computation}, where we want an efficient way to compute certain quantities like transition probabilities, or transition rates.

For \emph{fast exploration}, we herein propose the \emph{guided diffusion brigde} algorithm with the effective tracking choice 
\begin{equation}
  u(t,x) \;=\; J_\xi(x)^\top G_t\cdot\bigl(Z_t-\xi(x)\bigr),
  \label{eq:intro_tracking_control}
\end{equation} 
so that $\xi(X_t^{u})$ tracks a prescribed path
$Z_t$, e.g.,  given by the effective dynamics on $[T,T_{\mathrm{new}}]$, and $G_t\succeq 0$ is a (possibly time-dependent) gain (matrix or scalar).
Intuitively, \eqref{eq:intro_tracking_control} nudges the full system toward the CV path produced by the cheap effective dynamics, while the uncontrolled drift $b$ and the physical noise $\sigma$ still shape the
microscopic fluctuations. When $\xi$ is informative, the resulting guided trajectories can realize
physically plausible transition mechanisms in the full state space, even when the coarse path was
generated at much lower cost. Similar approaches utilizing guided diffusion systems have been studied in a general context, cf. \cite{Pedersen1995,DelyonHu2006,SchauerVanderMeulenVanZanten2017}, and for molecular dynamics \cite{DelarueKoehlOrland2016ConditionedLangevin,Orland2011LangevinBridges,DelarueKoehlOrland2017AbInitioCLD}. In these approaches, however, the diffusion bridge is not guided by a dynamically-informative effective dynamics but by endpoints or respective distributions. 

For \emph{accurate computation}, we will see that good choices of $\Psi$ in (\ref{eq:general_control}) have to be adapted to the specific dynamical quantities of interest, and that the literature on \emph{stochastic optimal control} (SOC) provides us with a strategy to choose $\Psi$ optimally. SOC allows to relate optimal control to variance-minimizing guidance, 
minimum-relative-entropy (Schr\"odinger bridge), and importance sampling in path space. This perspective has been discussed, e.g., in \cite{SchuetteKlusHartmann2023,Hartmannetal2026,HartmannKebirieial2019, HartmannSchuette2012OptimalNoneqForcing}.
It is discussed in more detail in Section~\ref{sec:soc-is} and clarifies the structure of optimal guidance and motivates future variance-reduction
strategies.

\paragraph{Correct-by-construction via Girsanov reweighting.}
Guidance introduces bias, so pathwise reweighting is required to recover statistics of the original
uncontrolled process conditioned on the coarse information. Under standard conditions, the
uncontrolled path measure on $[T,T_{\mathrm{new}}]$ is absolutely continuous with respect to the
controlled one, and the Radon--Nikodym derivative is given by the Girsanov formula
\begin{equation}
  \frac{d\mathbb{P}}{d\mathbb{P}^{u}}
  \;=\;
  \exp\!\Bigl(
   -\int_T^{T_{\mathrm{new}}} u(t,X_t^{u})\cdot dW_t
   -\int_T^{T_{\mathrm{new}}}\|u(t,X_t^{u})\|^2\,dt
  \Bigr).
  \label{eq:intro_girsanov}
\end{equation}
This yields weighted estimators for observables at $T_{\mathrm{new}}$ and, by resampling, a single
lifted state $X^\ast$. Importantly, the procedure is \emph{local}: it only requires evaluating $b(x)$,
$\xi(x)$, and $J_\xi(x)$ along the simulated guided paths, avoiding global sampling of fiber measures
(cf.\ the general emphasis on local Koopman/transfer-operator estimators in \cite{SchuetteKlusHartmann2023}).

\paragraph{Why this enables realistic full-system transition paths at reduced cost.}
The computational advantage comes from the separation of roles:
(i) the effective dynamics provides an inexpensive description of slow progress in CV space,
including coarse transition paths in the low dimensional latent space; (ii) the guided bridge uses this information to
steer the high-dimensional dynamics through the transition region while retaining microscopic
variability; and (iii) Girsanov reweighting and the related change of measure in path space corrects the steering bias, producing an importance-sampling
approximation of the desired conditional law. In the numerical experiments, this interplay is
illustrated by constructing full-dimensional transition paths whose CVs follow either a prescribed
coarse interpolation or a coarse trajectory generated by the effective model, while the resulting
full-system paths traverse the relevant transition regions in a dynamically consistent manner.

\paragraph{Contributions.}
The main contributions of this paper are:
\begin{enumerate}
\item a precise formulation of the lifting problem for coarse CV predictions in metastable diffusions;
\item a practical guided diffusion bridge algorithm that uses only local evaluations of $b$, $\xi$, and $J_\xi$
to generate full-dimensional path ensembles consistent with a given coarse reference path and allowing \emph{fast exploration};
\item an optimal choice of the guidance including a reweighting scheme enabling unbiased estimation of expectation values and transition quantities informed by the effective dynamics allowing for \emph{accurate computation};
\item numerical experiments demonstrating that, given an informative CV and a coarse effective trajectory,
the method produces realistic full-system transition paths and kinetic quantities at substantially reduced coarse-level effort.
\end{enumerate}

\paragraph{Outline.}
Section~\ref{sec:setting} introduces the setting and the available coarse information.
Section~\ref{sec:bidge} derives the guided bridge dynamics and the tracking control, summarizes the full lifting algorithm, discusses practical design choices and some options for  further improvements, while Section~\ref{sec:soc-is} focuses on deriving optimal guidance strategies for the accurate computation of specific dynamical quantities. 
Section~\ref{sec:numexp}  presents numerical experiments, and Section~\ref{sec:conclusion} concludes with a summary and some remarks on possible future extensions.

% % ------------------------------------------------------------------------

%%%%%%%%%%%%%%%%%%%%%%%%%%%%%%%%%%%%%%
%%%%%%%%%%%%%%%%%%%%%%%%%%%%%%%%%%%%%%
%%%%%%%%%%%%%%%%%%%%%%%%%%%%%%%%%%%%%%
%%%%%%%%%%%%%%%%%%%%%%%%%%%%%%%%%%%%%%

\section{Problem Setting and Available Coarse Information}\label{sec:setting}

We specify the full stochastic dynamics, the collective variable map $\xi$, and the associated coarse
(effective) dynamics on CV space. We then formalize the lifting task: given the last full state $X_T$ and a
coarse prediction in CV space over $[T,T_{\mathrm{new}}]$, construct full-dimensional
path ensembles that are consistent with the original dynamics and the prescribed coarse evolution.

\paragraph{Full-dimensional dynamics}
Let $X_t \in \mathcal{X}\subset \mathbb{R}^d$ follow an overdamped Langevin SDE
\begin{equation}\label{eq:fullSDE}
  dX_t = b(X_t)\,dt + \sigma\, dW_t,
  \qquad b(x) = -\nabla V(x),
\end{equation}
where $V:\mathcal{X}\to \mathbb{R}$ is smooth and confining, $\sigma>0$ is constant, and $W_t$ is a
$d$-dimensional Brownian motion. The slowest modes of the full dynamics, the rare transitions between its metastable sets, are given by the dominant eigenvalues and eigenfunctions of generator $\mathcal{L}$, see \cite{SchuetteKlusHartmann2023}, 
\begin{equation}\label{cL2}
\mathcal{L}f(x)=b(x)\cdot\nabla f(x)+\frac{\sigma^2}{2}\Delta f(x),
\end{equation}
or the associated transfer operator and there is a plethora of algorithms for computing these structures.

\paragraph{Collective variable and level sets}
Let $\xi:\mathcal{X}\to \mathbb{R}^m$ be a given smooth CV, with $m\ll d$, and define
\begin{equation}
  Z_t := \xi(X_t)\in \mathbb{R}^m, \qquad
  \mathbb{L}_z := \{x\in \mathcal{X}:\ \xi(x)=z\}.
\end{equation}
The level set $\mathbb{L}_z$ is a $(d-m)$-dimensional manifold under regularity of $\xi$.

\paragraph{Effective (coarse) dynamics}
An effective Markovian dynamics in CV space is assumed to be available,
\begin{equation}\label{eq:effSDE}
  dz_t = b^{\mathrm{e}}(z_t)\,dt + \hat{\sigma}(z_t)\,d\widetilde W_t,
\end{equation}
with $b^{\mathrm{e}}:\mathbb{R}^m\to\mathbb{R}^m$ and $\hat{\sigma}:\mathbb{R}^m\to\mathbb{R}^{m\times m}$ being determined, e.g., from conditional expectations along level sets (no global information required, as can be seen in \cite{eff-dyn-isolines}).  Sec.~\ref{app:ito} of the Appendix demonstrates how (\ref{eq:effSDE}) results from the full-dimensional dynamics.  We can advance \eqref{eq:effSDE} with a
coarse timestep $\Delta t$ much larger than the fine timestep $\delta t$ required for \eqref{eq:fullSDE}. Regarding notation, we will distinguish between $z_t$, denoting a trajectory of the effective dynamics (\ref{eq:effSDE}), and $Z_t=\xi(X_t)$, the representation of a full path $(X_t)$ in the latent space (as far as this distinction does not  disturb understanding).

\paragraph{Transition path theory}
The transition kinetics (dominant transition paths, rates, timescales) of the full system when starting in a set $A\subseteq\mathcal{X}$ and ending in another  disjoint set $B\subseteq \mathcal{X}$ can be characterized by means of transition path theory (TPT) \cite{E2005Transition,EVE2010TPTreview}. Its key concept, the \emph{committor function} $q(x)$ is the probability to reach $B$ before $A$ when starting from $x$:
\begin{equation}
q(x) = \mathbb{P}(\tau_x(B) < \tau_x(A)),\quad\text{for}\quad x\in \mathcal{X}\setminus(A\cup B),
\end{equation}
and can also be characterized as the unique solution of a boundary value problem of the full-system generator $\mathcal{L}$,
\begin{equation}\label{committor}
\mathcal{L} q(x) = 0 \;\text{for}\;  x\in \mathcal{X}\setminus(A\cup B),\quad\text{and}\quad 
q(x)=0, \;\text{for}\; x\in A,\quad
q(x)=1, \;\text{for}\; x\in B.
\end{equation}
TPT describes the ensemble of reactive trajectories— those paths that go from 
$A$ to $B$ without returning. For systems of form (\ref{eq:fullSDE}), the \emph{reactive density}
\begin{equation}\label{eq:reactiveDensity}
\mu_{AB}(x)=\mu(x)q(x)(1-q(x)),
\end{equation}
where $\mu(x)\propto \exp(-\beta V(x))$. Precisely, $\mu_{AB}$ is the stationary probability density of finding the system while on a reactive trajectory. The \emph{reactive probability flux} associated with transitions from $A$ to $B$ is
\begin{equation}\label{prob_flux}
j_{AB}(x)=\frac{1}{2}\mu(x)\,\sigma^2\,\nabla q(x),
\end{equation}
which allows to compute the associated transition rate and its main transition paths and channels, see \cite{MetznerSchuetteVandenEijnden2006TPTexamples} for illustrations and \cite{Sikorskietal2025} for a short overview including the adaptation of TPT to the effective dynamics.

% \paragraph{Information constraints}
% We assume:
% \begin{itemize}
% \item A previously computed full trajectory segment $\{X_t\}_{t\le T}$ is available, in particular the last state $X_T$.
% \item Starting from $Z_T=\xi(X_T)$, the effective dynamics has been advanced to time $T_{\mathrm{new}}=T+k\Delta t$,
% yielding coarse states $\{Z_{T+j\Delta t}\}_{j=0}^k$ and endpoint $Z^\ast:=Z_{T_{\mathrm{new}}}$.
% \item We do \emph{not} assume the ability to sample from the invariant measure or from the projected invariant measure in CV space.
% \end{itemize}
% The goal is to construct a random full state $X^\ast$ at time $T_{\mathrm{new}}$ such that $\mathcal{X}i(X^\ast)\approx Z^\ast$ and
% \[
%   \mathrm{Law}(X^\ast)\ \approx\ \mathrm{Law}\bigl(X_{T_{\mathrm{new}}}\ \big|\ \text{coarse information}\bigr).
% \]

%%%%%%%%%%%%%%%%%%%%%%%%%%%%%%%%%%%%%%
%%%%%%%%%%%%%%%%%%%%%%%%%%%%%%%%%%%%%%
%%%%%%%%%%%%%%%%%%%%%%%%%%%%%%%%%%%%%%
%%%%%%%%%%%%%%%%%%%%%%%%%%%%%%%%%%%%%%
%%%%%%%%%%%%%%%%%%%%%%%%%%%%%%%%%%%%%%
\section{Guided diffusion bridge method}\label{sec:bidge}

We now derive the controlled full-dimensional dynamics that will be used to lift coarse CV trajectories to
realistic transition paths of the original system. Starting from the baseline diffusion \eqref{eq:intro_full_sde},
we introduce a feedback control that steers the collective variable $\xi(X_t)$ toward a prescribed reference path $\bar z_t$ that was obtained from the effective dynamics, or denotes any other proposed path in latent space. We then discuss how this guidance is implemented in
discrete time and how the resulting bias is corrected later via Girsanov reweighting.

\subsection{Reference path in CV space}
Assume that from an effective simulation we have coarse states
\[
  z_T,\ z_{T+\Delta t},\ \dots,\ z_{T+k\Delta t}=z^\ast.
\]
We construct a continuous reference path $\bar z_t$ on $[T,T_{\mathrm{new}}]$ by interpolation. A simple
choice is piecewise linear interpolation: for $t\in[T+j\Delta t,\,T+(j+1)\Delta t]$,
\begin{equation}\label{eq:Zref}
  \bar z_t = z_{T+j\Delta t} + \frac{t-(T+j\Delta t)}{\Delta t}\,\bigl(z_{T+(j+1)\Delta t}-z_{T+j\Delta t}\bigr).
\end{equation}

\subsection{Controlled full dynamics (guidance)}
We introduce a controlled version of \eqref{eq:fullSDE} on $[T,T_{\mathrm{new}}]$:
\begin{equation}\label{eq:controlledSDE}
  dX^u_t = \bigl(b(X^u_t)+\sigma u(t,X^u_t)\bigr)\,dt + \sigma\,dW_t, \qquad X^u_T = X_T,
\end{equation}
where the guidance control $u$ is chosen so that the CV $\xi(X^u_t)$ tracks $\bar z_t$.
This tracking control is realized in the canonical way:
Let $J_\xi(x)\in \mathbb{R}^{m\times d}$ denote the Jacobian of $\xi$, that is,
\begin{equation}\label{eq:J_xi}
J_{\xi}(x)
:= \left[\frac{\partial \xi_i(x)}{\partial x_j}\right]_{i=1,\ldots,m}^{j=1,\ldots,d}
\in \mathbb{R}^{m\times d}.
\end{equation}
and set the control to
\begin{equation}\label{eq:control}
  u(t,x) = J_\xi(x)^\top G_t\,\bigl(\bar z_t-\xi(x)\bigr),
\end{equation}
where $G_t\in \mathbb{R}^{m\times m}$ is symmetric positive definite (e.g.\ $G_t=\kappa I_m$).

% \paragraph{Ito linearization.}
% Applying It\^{o}'s formula to $Z^u_t:=\xi(X^u_t)$ gives
% \begin{equation}\label{eq:itoZ}
%   dZ^u_t
%   = J_\xi(X^u_t)\bigl(b(X^u_t)+u(t,X^u_t)\bigr)\,dt
%     + \frac{\sigma^2}{2}\,\Delta \xi(X^u_t)\,dt
%     + \sigma\,J_\xi(X^u_t)\,dW_t.
% \end{equation}
% Substituting \eqref{eq:control} yields a stabilizing drift term pulling $\xi(X^u_t)$ towards $\bar Z_t$.

Numerical simulation of (\ref{eq:controlledSDE}) can utilize standard Euler--Maruyama discretization in time:
Choose a fine timestep $\delta t$ with $\delta t\ll \Delta t$ and set $M=(T_{\mathrm{new}}-T)/\delta t$.
Let $t_n=T+n\delta t$. For each trajectory we simulate
\begin{equation}\label{eq:EM}
  X^u_{n+1}
  = X^u_n + \Bigl(b(X^u_n)+\sigma u(t_n,X^u_n)\Bigr)\delta t + \sigma\sqrt{\delta t}\,\eta_n,
\end{equation}
with i.i.d.\ $\eta_n\sim \mathcal{N}(0,I_d)$.

\subsection{Girsanov reweighting}
Under standard conditions (e.g.\ Novikov), the path measure of the uncontrolled
process on $[T,T_{\mathrm{new}}]$ is absolutely continuous with respect to the controlled one, with
Radon--Nikodym derivative
\begin{equation}\label{eq:girsanov}
  \frac{d\mathbb{P}}{d\mathbb{P}^u}
  = \exp\!\left(
    - \int_T^{T_{\mathrm{new}}} u(t,X^u_t)\cdot dW_t
    - \int_T^{T_{\mathrm{new}}} \|u(t,X^u_t)\|^2\,dt
  \right).
\end{equation}

\paragraph{Discrete weight formula.}
For a simulated guided trajectory $\{X^u_n\}_{n=0}^M$, define increments $\Delta W_n=\sqrt{\delta t}\,\eta_n$ and
control values $u_n=u(t_n,X^u_n)$. A standard discretization of \eqref{eq:girsanov} is
\begin{equation}\label{eq:discweights}
  w \approx \exp\!\left(
    - \sum_{n=0}^{M-1} u_n\cdot \Delta W_n
    - \sum_{n=0}^{M-1} \|u_n\|^2\,\delta t
  \right).
\end{equation}

\subsection{Estimators and lifted samples}
Let $f:\mathcal{X}\to\mathbb{R}$ be an observable. If $\{X^{u,(j)}_M\}_{j=1}^N$ are endpoints of $N$ guided
trajectories with weights $\{w^{(j)}\}_{j=1}^N$, the importance-sampling estimator is
\begin{equation}\label{eq:IS}
  \mathbb{E}\bigl[f(X_{T_{\mathrm{new}}})\ \big|\ \text{coarse constraint}\bigr]
  \approx
  \frac{\sum_{j=1}^N w^{(j)} f\!\left(X^{u,(j)}_M\right)}{\sum_{j=1}^N w^{(j)}}.
\end{equation}
To output a single lifted full state, resample an index $\ell$ from normalized weights
$\tilde w^{(j)}=w^{(j)}/\sum_\ell w^{(\ell)}$ and set $X^\ast=X^{u,(\ell)}_M$.

\subsection{Algorithm}

Before stating the full scheme, we briefly summarize the lifting procedure carried out in
Algorithm~1. Assume we are given the last available full state $X_T$ at time $T$ and, on the
time interval $[T,T_{\mathrm{new}}]$, a coarse trajectory in CV space
$\{Z_{T+j\Delta t}\}_{j=0}^k$ (e.g.\ produced by an effective dynamics), where $T_{\mathrm{new}}=T+k\Delta t$.
We first construct a continuous reference path $\bar z_t$ on $[T,T_{\mathrm{new}}]$ by interpolating the
coarse points. We then generate an ensemble of $N$ full-dimensional paths by simulating a
\emph{guided} version of the full dynamics: at each fine timestep $\delta t$ we add a feedback control of the form described above
which steers the CV $\xi(X_t)$ toward the reference $\bar z_t$ while preserving the stochasticity of the
full model. Since the guidance biases the path distribution, we compute for each simulated trajectory
a Girsanov log-weight by accumulating its discrete-time analogue along the path. Finally, the weights
are normalized and used either (i) to form weighted estimators of observables at $T_{\mathrm{new}}$, or
(ii) to resample a representative lifted endpoint $X^\ast$ (and, if desired, a representative lifted path).
Importantly, the method is local and on-the-fly: it requires only evaluations of $b(x)$, $\xi(x)$, and
$J_\xi(x)$ along the guided trajectories and avoids global sampling of conditional equilibrium measures
on CV level sets.

\medskip
Algorithm~1 below makes these steps explicit: (1) interpolate the coarse CV path to obtain $\bar z_t$;
(2) simulate $N$ guided Euler--Maruyama trajectories of the controlled full dynamics; (3) update the
Girsanov log-weights in each step; and (4) normalize weights and resample to obtain a lifted endpoint.

\begin{algorithm}[h]
\caption{Guided Langevin Bridge Lifting (Effective-Path Guidance + Girsanov Reweighting)}
\begin{algorithmic}[1]
\Require Full drift $b(x)=-\nabla V(x)$, noise $\sigma$, CV $\xi$ with Jacobian $J_\xi$.
\Require Last full state $X_T$, effective states $\{Z_{T+j\Delta t}\}_{j=0}^k$,
        $T$, $T_{\mathrm{new}}=T+k\Delta t$.
\Require Fine timestep $\delta t$ (with $\delta t\ll\Delta t$), number of guided trajectories $N$,
        control gain schedule $G_t$.
\Ensure Weighted ensemble $\{(X^{u,(j)}_{T_{\mathrm{new}}},w^{(j)})\}_{j=1}^N$ and/or a lifted sample $X^\ast$.
\State Construct reference path $\bar z_t$ on $[T,T_{\mathrm{new}}]$ by interpolation of $\{Z_{T+j\Delta t}\}$.
\For{$j=1$ to $N$}
  \State Set $X^{u,(j)}_0 \gets X_T$, $\log w^{(j)} \gets 0$.
  \For{$n=0$ to $M-1$ where $M=(T_{\mathrm{new}}-T)/\delta t$}
    \State $t_n \gets T+n\delta t$.
    \State $u_n \gets J_\xi(X^{u,(j)}_n)^\top G_{t_n}\bigl(\bar Z_{t_n}-\xi(X^{u,(j)}_n)\bigr)$.
    \State Sample $\eta_n \sim \mathcal{N}(0,I_d)$ and set $\Delta W_n \gets \sqrt{\delta t}\,\eta_n$.
    \State $X^{u,(j)}_{n+1} \gets X^{u,(j)}_n + (b(X^{u,(j)}_n)+u_n)\delta t + \sigma \Delta W_n$.
    \State $\log w^{(j)} \gets \log w^{(j)} - \frac{1}{\sigma} u_n\cdot \Delta W_n - \frac{1}{2\sigma^2}\|u_n\|^2\delta t$.
  \EndFor
  \State $X^{u,(j)}_{T_{\mathrm{new}}} \gets X^{u,(j)}_M$, \quad $w^{(j)}\gets \exp(\log w^{(j)})$.
\EndFor
\State Normalize weights $\tilde w^{(j)} \gets w^{(j)}/\sum_{\ell=1}^N w^{(\ell)}$.
\State \textbf{Option A (ensemble):} return $\{(X^{u,(j)}_{T_{\mathrm{new}}},\tilde w^{(j)})\}_{j=1}^N$.
\State \textbf{Option B (single lift):} sample $J$ with $\mathbb{P}(J=j)=\tilde w^{(j)}$ and set $X^\ast\gets X^{u,(J)}_{T_{\mathrm{new}}}$.
\end{algorithmic}
\end{algorithm}

%\subsection{Design choices and practical refinements}

\paragraph{Choosing the gain $G_t$}
Large gains enforce $\xi(X^u_t)\approx \bar z_t$ strongly but can lead to weight degeneracy.
Small gains yield better weight balance but weaker tracking.
A diagnostic is the effective sample size
\[
\mathrm{ESS} = \frac{1}{\sum_{j=1}^N (\tilde w^{(j)})^2}.
\]
If ESS collapses, reduce the gain or increase $N$. Gain schedules that increase toward $T_{\mathrm{new}}$
often improve endpoint accuracy without excessive early forcing.

\paragraph{Numerical stability}
Ill-conditioning of $J_\xi J_\xi^\top$ can lead to stiff feedback. Remedies include:
\begin{itemize}
\item Preconditioning: replace $G_t$ by $G_t (J_\xi J_\xi^\top + \lambda I_m)^{-1}$, $\lambda>0$.
\item Clipping: bound $\|u(t,x)\|$ to avoid extreme controls.
\item Smaller $\delta t$ when the control becomes stiff.
\end{itemize}

\paragraph{Sequential Monte Carlo variant}
To combat weight degeneracy over long horizons, incorporate resampling and mutation steps (SMC):
propagate guided particles for several steps, compute incremental weights, resample if ESS drops below a
threshold, then continue. This yields a particle approximation of the conditioned path measure and is
often more robust for large $T_{\mathrm{new}}-T$.

%%%%%%%%%%%%%%%%%%%%%%%%%%%%%%%%%%%%%%%%%%
%%%%%%%%%%%%%%%%%%%%%%%%%%%%%%%%%%%%%%%%%%
\section{Optimal Guidance}
\label{sec:soc-is}

Next, we want to discuss the question of how to choose the guidance optimally. To this end, we briefly summarize the so-called \emph{stochastic optimal control} (SOC) viewpoint that motivates drift-guided sampling in the guided-bridge construction. In Sec.~\ref{ssec:soc} we review recent literature, see for \cite{SchuetteKlusHartmann2023,Hartmannetal2026,HartmannKebirieial2019, HartmannSchuette2012OptimalNoneqForcing} for example, where the SOC perspective, its equivalent change of measure formulation, and the related theory of zero-variance estimators and optimal importance sampling are discussed. In Sec.~\ref{subsec:bk-hjb-effective}, we then consider the SOC for the effective dynamics, and utilize the results to present the optimal guidance scheme in Sec.~\ref{ssec:optGuidance}. 

\subsection{SOC formulation and optimal forcing}\label{ssec:soc}
Let $(X_s)_{s\in[t,T]}$ follow the uncontrolled diffusion
\begin{equation}
  dX_s = b(X_s)\,ds + \sigma(X_s)\,dW_s,\qquad X_t=x,
  \label{eq:soc-uncontrolled}
\end{equation}
with state-dependent diffusion $\sigma=\sigma(x)$, and consider controlled dynamics obtained by adding a progressively measurable forcing $u_s\in\mathbb{R}^m$ through the diffusion channels,
\begin{equation}
  dX_s^{u} = \bigl(b(X_s^{u})+\sigma(X_s^{u})u_s\bigr)\,ds + \sigma(X_s^{u})\,dW_s,
  \qquad X_t^{u}=x.
  \label{eq:soc-controlled}
\end{equation}
We assume that our aim is to compute the following (path-dependent) observable: 
\begin{equation}
  \Psi(t,x)
  \;:=\;
  -\log \mathbb{E}_{t,x}\!\left[\exp\!\bigl(-g(X_T)\bigr)\right],
  \label{eq:Psi-def}
\end{equation}
where $g:\mathbb{R}^d\to\mathbb{R}\cup\{+\infty\}$ is called the \emph{terminal cost}. Interestingly, there is a standard variational (log-transform) identity that yields that the observable $\Psi$ admits an equivalent optimal control representation
\begin{equation} \label{eq:SOC}
  \Psi(t,x)
  \;=\;
  \inf_{u}\;
  \mathbb{E}^{u}_{t,x}\!\left[
    \frac12\int_t^T \|u_s\|^2\,ds
    \;+\;
    g\!\bigl(X_T^{u}\bigr)
  \right],
\end{equation}
where $\mathbb{E}^{u}_{t,x}$ denotes expectation under the guided/controlled diffusion \eqref{eq:soc-controlled}. The observable $\Psi$ is therefore also called the \emph{SOC functional}. The optimizer $u^\ast$ in (\ref{eq:SOC}) is the \emph{optimal forcing}; it trades off steering the dynamics toward low terminal cost against the quadratic control energy.

\paragraph{Explicit formula for the optimal control.}
 In the control setting \eqref{eq:soc-controlled}, the optimal forcing admits the explicit gradient form
\begin{equation}
  u^\ast(t,x) \;=\; -\,\sigma(x)^\top \nabla_x \Psi(t,x).
  \label{eq:u-star}
\end{equation}
This is the control-theoretic target that the guided-bridge forcing in Algorithm~1 aims to approximate in a low-dimensional, computationally tractable way (via the chosen reaction coordinate $\xi$ and a gain schedule), see below.

\paragraph{Change of measure and importance sampling}
The drift modification in \eqref{eq:soc-controlled} corresponds to a change of path measure (Girsanov, cf. (\ref{eq:girsanov})). Writing $P$ for the law of \eqref{eq:soc-uncontrolled} and $P^{u}$ for the law of \eqref{eq:soc-controlled} on $[t,T]$, the Radon--Nikodym derivative satisfies
\begin{equation}
  \frac{dP}{dP^{u}}
  \;=\;
  \exp\!\left(
    -\int_t^T u_s\cdot dW_s
    \;-\;
    \frac12\int_t^T \|u_s\|^2\,ds
  \right).
  \label{eq:girsanov}
\end{equation}
Consequently, any expectation under $P$ can be evaluated under $P^{u}$ with reweighting. In particular, for the terminal weight $\exp(-g(X_T))$,
\begin{equation}
  \mathbb{E}_{t,x}\!\left[e^{-g(X_T)}\right]
  \;=\;
  \mathbb{E}^{u}_{t,x}\!\left[
    e^{-g(X_T^{u})}\,
    \exp\!\left(
      -\int_t^T u_s\cdot dW_s
      \;-\;
      \frac12\int_t^T \|u_s\|^2\,ds
    \right)
  \right].
  \label{eq:is-identity}
\end{equation}
Thus, choosing $u$ is simultaneously (i) choosing a \emph{control} in \eqref{eq:soc-controlled} and (ii) choosing an \emph{importance-sampling measure} in \eqref{eq:is-identity}. The optimizer $u^\ast$ (formally) produces the variance-minimizing change of measure within the class of drift changes $b\mapsto b+\sigma u$.

\paragraph{Example (Transition Probability):}
Let $B\subset\mathbb{R}^d$ be measurable and define the finite-time probability
\begin{equation}
  p_B(t,x) \;:=\; \mathbb{P}(X_T\in B\mid X_t=x)=\mathbb{E}_{t,x}\!\left[\mathbf{1}_B(X_T)\right].
  \label{eq:pB}
\end{equation}
This fits into \eqref{eq:Psi-def} by choosing the (hard) terminal cost
\begin{equation}
  g(x)= -\log \mathbf{1}_B(x)
  \;=\;
  \begin{cases}
    0, & x\in B,\\
    +\infty, & x\notin B,
  \end{cases}
\end{equation}
(or a smooth penalty approximation), where $\mathbf{1}_B$ denote the indicator function of set $B$. Then $\exp(-g(X_T))=\mathbf{1}_B(X_T)$ and
\begin{equation}
  \Psi(t,x)
  \;=\;
  -\log p_B(t,x)
  \;=\;
  \inf_{u}\;
  \mathbb{E}^{u}_{t,x}\!\left[
    \frac12\int_t^T \|u_s\|^2\,ds
    \;-\;
    \log \mathbf{1}_B\!\bigl(X_T^{u}\bigr)
  \right].
  \label{eq:Psi-event}
\end{equation}
In this way, estimating $p_B(t,x)$ becomes equivalent to a finite-time SOC problem whose optimal forcing $u^\ast$ makes the event $X_T\in B$ typical under $P^{u^\ast}$. (

%%%%%%%%%%%%%%%%%%%%%%%%%%%%%
%%%%%%%%%%%%%%%%%%%%%%%%%%%%%%%%%%%%%%%%%
%%%%%%%%%%%%%%%%%%%%%%%%%%%%%%%%%%%%%%%%%
\paragraph{Random stopping times:}
If the terminal time is replaced by a stopping time $\tau$ (e.g.\ a first exit/hitting time), it is natural to consider payoffs that include a \emph{path integral} up to $\tau$. A standard choice is the Feynman--Kac type functional
\begin{equation}
  \Phi(x)
  \;:=\;
  \mathbb{E}_x\!\left[
    \exp\!\Bigl(-\int_{0}^{\tau} f(X_s)\,ds\Bigr)\,e^{-g(X_\tau)}
  \right],
  \label{eq:fk-stopping}
\end{equation}
where $f\ge 0$ plays the role of a running cost (or killing/discount rate) and $e^{-g}$ is a boundary payoff.

If $\tau$ denotes the first exit time from a domain $D$, i.e.\ $\tau=\inf\{t\ge 0: X_t\notin D\}$.
Then the function $\Phi$ in \eqref{eq:fk-stopping} solves the boundary value problem
\begin{equation}
  (\mathcal{L}\Phi)(x) \;-\; f(x)\,\Phi(x)\;=\;0 \quad \text{in } D,
  \qquad
  \Phi(x)=e^{-g(x)}\quad \text{on } \partial D,
  \label{eq:bk-fk-stopping}
\end{equation}
where $\mathcal{L}$ denotes the generator of the uncontrolled diffusion. Defining the log-value
\[
\Psi(x):=-\log \Phi(x),
\]
one obtains the random-horizon SOC representation
\begin{equation}
  \Psi(x)
  \;=\;
  \inf_{u}\;
  \mathbb{E}_x^{u}\!\left[
    \int_0^{\tau^u}\Bigl(f(X_s^{u})+\tfrac12\|u_s\|^2\Bigr)\,ds
    \;+\;
    g(X_{\tau^u}^{u})
  \right],
  \label{eq:soc-stopping-fk}
\end{equation}
with controlled dynamics given by (\ref{eq:soc-controlled}) and its hitting time $\tau^u$, and the optimal control again has the from $u^\ast(x)=-\sigma(x)^\top\nabla\Psi(x)$.

% \paragraph{Remark:}
% In addition to the exponential (Feynman--Kac) functional \eqref{eq:fk-stopping}, one often considers the \emph{additive} random-horizon payoff
% \begin{equation}
%   \Phi(x)
%   \;:=\;
%   \mathbb{E}_x\!\left[
%     \int_{0}^{\tau} f(X_s)\,ds \;+\; G(X_\tau)
%   \right],
%   \label{eq:additive-stopping}
% \end{equation}
% where $f$ is a running cost accumulated along the path up to $\tau$, and $G$ prescribes a boundary payoff at the stopping time. If $\tau$ is the first exit time from a domain $D$, then $\Phi$ solves the (stationary) backward equation
% \begin{equation}
%   (\mathcal{L}\Phi)(x) \;=\; -\,f(x)\quad \text{in } D,
%   \qquad
%   \Phi(x)=G(x)\quad \text{on } \partial D,
%   \label{eq:poisson-stopping}
% \end{equation}
% i.e.\ a Poisson problem with Dirichlet boundary data on $\partial D$. 

%%%%
\paragraph{Example (Committor function):}
The committor
\[
q(x):=\mathbb{P}_x(\tau_B<\tau_A),
\]
with $A,B$ being disjoint sets with stopping times $\tau_A,\tau_B$, is the solution of (\ref{committor}). Let $\tau=\tau_{A\cup B}=\min(\tau_A,\tau_B)$ denote the stopping time of $A\cup B$ and $\partial B$ the (smooth) boundary of $B$ (where hitting $B$ takes place), then
\[
q(x)=\mathbb{E}_x\!\left[
    \exp\!\Bigl(-\int_{0}^{\tau} f(X_s)\,ds\Bigr)\,e^{-g(X_\tau)},
  \right],
\]
with $f=0$ and $g=-\log\mathbf{1}_{\partial B}$ and consequently
\begin{equation}\label{eq:committor_optimal_forcing}
-\log q(x) = \;=\;
  \inf_{u}\;
  \mathbb{E}_x^{u}\!\left[\tfrac12
    \int_0^{\tau^u}\|u_s\|^2\,ds
    \;-\;
    \log \mathbf{1}_{\partial B}(X_{\tau^u}^{u})
  \right],
\end{equation}
with optimal force 
\[
u^\ast(x)=\sigma(x)^\top \nabla \log q(x).
\]

%%%%%%%%%%%%%%%%%%%%%%%%%%%%%%%%%%%%%%%%%
%%%%%%%%%%%%%%%%%%%%%%%%%%%%%%%%%%%%%%%%%

%%%%%%%%%%%%%%%%%%%%%%%%%%%%%%%%%%%%%%%%%%%%

\subsection{Effective dynamics}
\label{subsec:bk-hjb-effective}

Consider the effective one-dimensional diffusion (like (\ref{eq:effSDE}) with $D(z)=\frac12 \hat{\sigma}^2$ in 1d),
\begin{equation}
  dZ_s \;=\; b(Z_s)\,ds \;+\; \sqrt{2D(Z_s)}\,dW_s,\qquad Z_0=z_0,
\end{equation}
with drift $b(z)=c+\lambda z$ and (state-dependent) diffusion coefficient $D(z)=D_{\mathrm{eff}}(z)$. For a threshold $z_\ast\in(0,1)$ define the target set
\(
B_z=\{z\in[0,1]: z>z_\ast\}
\)
and the finite-time probability
\begin{equation}\label{p(t,z)}
  p(s,z)\;:=\;\mathbb{P}\bigl(z_t\in B_z \,\big|\, z_s=z\bigr)
  \;=\;\mathbb{E}\bigl[\mathbf{1}_{\{z_t>z_\ast\}}\mid z_s=z\bigr],
  \qquad 0\le s\le t.
\end{equation}
Then $p$ satisfies the backward Kolmogorov equation
\begin{equation}
  \partial_s p(s,z)\;+\; \mathcal{L}_{\text{eff}}\, p(s,z)\;=\;0,
  \qquad
  p(t,z)=\mathbf{1}_{\{z>z_\ast\}},
  \label{eq:bk-effective}
\end{equation}
supplemented with the boundary condition consistent with the effective dynamics on $[0,1]$ (no-flux), and using the generator of the effective dynamics,
\[
\mathcal{L}_{\text{eff}} = b(z)\,\partial_z 
  \;+\; D(z)\,\partial_{zz}^2.
\]
To cast this into an SOC form, introduce the value function
\begin{equation}
  \Psi(s,z)\;:=\;-\log p(s,z),
  \qquad\text{so that}\qquad
  p(s,z)=e^{-\Psi(s,z)}.
\end{equation}
Using the controlled effective dynamics (control through the diffusion channels)
\begin{equation}
  dZ_s^{u} \;=\; b(Z_s^{u})\,ds \;+\; \sqrt{2D(Z_s^{u})}\,\bigl(u_s\,ds+dW_s\bigr),
  \label{eq:controlled-effective}
\end{equation}
the SOC functional for the event $\{z_t>z_\ast\}$ is
\begin{equation}
  \Psi(s,z)\;=\;\inf_{u}\;
  \mathbb{E}^{u}_{s,z}\!\left[
    \frac12\int_s^{t} u_r^2\,dr \;+\; g(z_t^{u})
  \right],
  \qquad
  g(z)= -\log \mathbf{1}_{\{z>z_\ast\}},
  \label{eq:Psi-effective-event}
\end{equation}
(with $g$ interpreted as a hard constraint, or replaced by a smooth penalty approximation in computations).

% The corresponding Hamilton--Jacobi--Bellman equation for $\Psi$ reads
% \begin{equation}
%   \partial_s \Psi(s,z)
%   \;+\; \mathcal{L}_{\text{eff}} \Psi(s,z)
%   \;-\; D(z)\,\bigl(\partial_z \Psi(s,z)\bigr)^2
%   \;=\; 0.
%   \qquad
%   \Psi(t,z)=g(z).
%   \label{eq:hjb-effective}
% \end{equation}
The optimal feedback control of the effective dynamics is explicit:
\begin{equation}
  \ueff(s,z)\;=\;-\sqrt{2D(z)}\,\partial_z \Psi(s,z)
  \;=\;\sqrt{2D(z)}\,\frac{\partial_z p(s,z)}{p(s,z)}.
  \label{eq:optimal-control-effective}
\end{equation}
In the present work, the guided-bridge forcing used for lifting can be understood as a structured approximation of such an optimal feedback in terms of the chosen reaction coordinate and a time-dependent gain schedule.

It should be obvious that the formula (\ref{eq:committor_optimal_forcing}) for using optimal control in the case of the committor function can be transferred to the effective dynamics.

%%%%%%%%%%%%%%%%%%%%%%%%%%%%%%%%%%%%%%%%%%%%%%%%%%%%%%%%%%%%
%%%%%%%%%%%%%%%%%%%%%%%%%%%%%%%%%%%%%%%%%%%%%%%%%%%%%%%%%%%%
\subsection{Optimal guidance for the full-dimensional system}\label{ssec:optGuidance}
Returning to the controlled full-dimensional system (\ref{eq:soc-controlled}), we may choose a control $\tu$ that, since the computation of the observable $\Psi$ is infeasible or too expensive, is based on a related observable $\tilde{\Psi}\circ \xi$ instead, where $\xi$ is our dynamically-informative CV. The observable $\tilde{\Psi}$ thus depends on the effective variable $z=\xi(x)$ alone. This leads to a control
of the form
\begin{equation}\label{eq:optimal_guidance}
    \tu(s,x) = -\sigma^\top (\nabla_x\xi)\; \partial_z\tilde{\Psi}(s,z) =\frac{\sigma^\top}{\sqrt{2D(z)}} J_\xi(x)^\top \ueff(s,z),
\end{equation}
with $z=\xi(x)$ and $J_\xi$ denoting the Jacobian of $\xi$ as introduced in (\ref{eq:J_xi}).
In general, $\tu$ is not identical to the optimal control $u^\ast$  of the full system as given by (\ref{eq:u-star}) for the observable $\Psi$, but an \emph{optimal guidance} based on the optimal control $\ueff$ of the effective dynamics as given by (\ref{eq:optimal-control-effective}) for the respective observable $\tilde{\Psi}$. The resulting \emph{optimally guided diffusion} reads
\[
  dX_s^{u} = \bigl(b(X_s^{u})-\sigma(X_s^{u})\sigma(X_s^{u})^\top \;\partial_z\tilde{\Psi}(\xi(X_s))\,J_\xi(X_s^{u})^\top\bigr)\,ds + \sigma(X_s^{u})\,dW_s,
  \qquad X_0^{u}=x.
\]
A similar idea has already been pursued in \cite{HartmannSchuetteZhang2016Nonlinearity} where it was used to construct full-dimensional ansatz functions adapted to coarse information; herein, we go one step further by avoiding to compute the optimal force for the full-dimensional system instead of using the coarse information for \emph{guiding} it.

\paragraph{Transition Probability:} Let $\xi$ be 1-d. Then, for the finite-time probability and $\Psi(s,x)=-\log p_B(s,x)$ with $B=\{x:\, \xi(x)>z_\ast\}$, for example, we may choose
\[
\tilde{\Psi}(s,z)=-\log p(s,z)
\]
with $p(s,z)$ as defined in (\ref{p(t,z)}).
Returning to the case of a scalar constant noise intensity $\sigma(x)=\sigma$,  the optimally guided diffusion for computing $p_B(t,x_0)$ reads
\begin{equation}
  dX_s^{u} = \bigl(b(X_s^{u})+\sigma^2 \Big(\partial_z \log p(s,\xi(X^{u}_s))\Big)\, J_\xi(X_s^{u})^\top\bigr)\,ds + \sigma\,dW_s,
  \qquad X_0^{u}=x,
  \label{eq:optGuided}
\end{equation}
where we can get $p(s,z)$ by solving the PDE (\ref{eq:bk-effective}) which lives in 1-d if $\xi$ is a one-dimensional CV. The resulting estimate for the desired quantity is
\begin{equation}\label{eq:MC_estimate_pB}
\hat{p}_B(T,x) = \exp\left(-\mathbb{E}^{\tu}_{0,x}\left[ \frac12 \int_0^T \|\tu_s\|^2 ds -  \log \mathbf{1}_B(X^{\tu}_T)\right]\right),
\end{equation}
where 
\[
\int_0^T \|\tu_s\|^2 ds = \sigma^2 \int_0^T \|J_\xi(X_s^{u})\|^2 \Big(\partial_z \log p(s,\xi(X^{u}_s))\Big)^2 ds,
\]
where $p(s,z)$ is given via the backward Kolmogorov PDE (\ref{eq:bk-effective}) of the effective dynamics, or, respectively, via 
\[
p(s,z)=\exp((t-s)\cL_{\text{eff}})\mathbf{1}_{z>z_\ast}=\mathcal{K}^{t-s}_{\text{eff}}\mathbf{1}_{z>z_\ast},
\]
where $\mathcal{K}^t_{\text{eff}}=\exp(t\cL_{\text{eff}})$ denotes the Koopman or transfer operator of the effective dynamics.

The control $\tu$ used in the optimal guidance scheme is not the optimal control $u^\ast$ of (\ref{eq:u-star}) but an approximation of it. Therefore, it may happen that some of the guided paths will not reach the target set.   For practical reasons, we thus introduce a boosting factor $\kappa$ and replace $\tilde{u}$ by $\kappa\tilde{u}$, as well as the regularized estimator $\hat{p}_{B,\epsilon}(T,x)=\max(\hat{p}_B(T,x),\epsilon)$ with small $\epsilon>0$. 

%%%%%%%%%%%%%%%%%%%%
\paragraph{Committor:} 
%%%%%%%%%%%%%%%%%%%%
If the committor $q$ of the full system is known or can be approximated in terms of a good collective variable $\xi$ as $q(x)=\tilde q(\xi(x))$, then
\[
\nabla\log q(x)=\frac{\tilde q'(\xi(x))}{\tilde q(\xi(x))}\,\nabla\xi(x)=J_\xi(x)^\top \partial_z \log \tilde q(\xi(x)),
\]
and the guiding control (\ref{eq:optimal_guidance}) admits the explicit CV-based form
\begin{equation}\label{eq:committor_guidance}
\tu(x)
=\sigma(x)^\top J_\xi(x)^\top\;\frac{\tilde q'(\xi(x))}{\tilde q(\xi(x))},
\end{equation}
which leads us to the committor estimate
\begin{equation}\label{eq:committor_controlled}
\hat{q}(x) = \exp\Big(-\mathbb{E}_x^{\tu}\!\left[\tfrac12
    \int_0^{\tau^{\tu}}\|\tu_s\|^2\,ds
    \;-\;
    \log \mathbf{1}_{\partial B}(X_{\tau^{\tu}}^{u})
  \right]\Big).
\end{equation}

%%%%%%%%%%%%%%%%%%%%%%%%%%%%%%%%%%%%%%%%%%%%%%
%%%%%%%%%%%%%%%%%%%%%%%%%%%%%%%%%%%%%%%%%%%%%%%

\subsection{Policy iteration}
The value function $\Psi$ solves a nonlinear Hamilton--Jacobi--Bellman (HJB) equation. With $a(x)=\sigma(x)\sigma(x)^\top$ and generator $\mathcal{L}f=b\cdot\nabla f + \tfrac12 a:\nabla^2 f$, the corresponding HJB equation reads 
\begin{equation}
  \partial_t \Psi(t,x)
  \;+\;
  \mathcal{L}\Psi(t,x)
  \;-\;
  \frac12\bigl\|\sigma(x)^\top \nabla_x \Psi(t,x)\bigr\|^2
  \;=\; 0,
  \qquad
  \Psi(T,x)=g(x).
  \label{eq:HJB}
\end{equation}
Policy iteration (Howard's algorithm) computes $\Psi$ and $u^\ast$ by alternating:
\begin{enumerate}
\item \textbf{Policy evaluation.} Given a current feedback $u_k(t,x)$, solve the \emph{linear} PDE
\begin{equation}
  \partial_t \Psi_k
  \;+\;
  \bigl(b+\sigma u_k\bigr)\cdot\nabla \Psi_k
  \;+\;
  \tfrac12 a:\nabla^2 \Psi_k
  \;+\;
  \frac12\|u_k\|^2
  \;=\; 0,
  \qquad
  \Psi_k(T,\cdot)=g(\cdot).
  \label{eq:policy-eval}
\end{equation}
\item \textbf{Policy improvement.} Update the control by pointwise minimization of the Hamiltonian, which in the quadratic case yields the explicit rule
\begin{equation}
  u_{k+1}(t,x)
  \;=\;
  -\,\sigma(x)^\top \nabla_x \Psi_k(t,x).
  \label{eq:policy-improve}
\end{equation}
\end{enumerate}
Under suitable conditions, $u_k\to u^\ast$ and $\Psi_k\to \Psi$. In high-dimensional molecular systems one typically replaces \eqref{eq:policy-eval} by an approximate evaluation (e.g.\ Monte Carlo regression under the controlled dynamics), while retaining the explicit improvement step \eqref{eq:policy-improve}; this yields an \emph{approximate policy iteration} scheme. In our setting, the guided-bridge forcing used in Algorithm~1 can be interpreted as a structured, low-dimensional surrogate for the optimal feedback \eqref{eq:u-star}, built from the reaction coordinate $\xi$ and a gain schedule. One of the vaiants closest to our approach is the \emph{Adaptive Biasing Force (ABF)} approach. ABF aims to compute the potential of mean force (free energy) related to the ABF $\xi$,
\[
A(z) \;=\; -\beta^{-1}\log\!\int e^{-\beta V(x)}\,\delta(\xi(x)-z)\,dx,
\qquad \beta=(k_BT)^{-1},
\]
by estimating its derivative (the mean force) along level sets of $\xi$.
In practice one constructs an on-the-fly estimator $\widehat{A}'(z)\approx A'(z)$ from time averages of suitable instantaneous force estimators conditioned on $\xi(x)\approx z$, and applies a compensating biasing force
\[
F_{\mathrm{bias}}(x)\;=\;-\widehat{A}'(\xi(x))\,\nabla \xi(x),
\]
so that the \emph{average} drift along $\xi$ is progressively cancelled and barriers in $A$ are flattened.
In the long-time limit, the resulting biased dynamics yields (approximately) uniform sampling in $z$ and a consistent reconstruction of $A$ via integration of $\widehat{A}'$. \cite{DarvePohorille2001,ComerEtAl2015ABF}

% --- Short overview over methods (as mentioned earlier) ---
\paragraph{Overview of common enhanced sampling / rare-event schemes.}
Umbrella sampling enforces exploration of prescribed CV windows and is typically combined with WHAM for unbiased free-energy reconstruction \cite{TorrieValleau1977,KumarEtAl1992WHAM}.
Metadynamics adds a history-dependent bias to escape metastable basins and recover free-energy profiles \cite{LaioParrinello2002}.
ABF estimates and cancels the mean force along a CV to obtain the potential of mean force \cite{DarvePohorille2001,ComerEtAl2015ABF}.
Variationally Enhanced Sampling (VES) selects a bias by minimizing a variational functional within a parametrized family \cite{ValssonParrinello2014VES}.
OPES reframes adaptive biasing in terms of targeting probability distributions, yielding robust on-the-fly bias construction \cite{InvernizziParrinello2020OPES}.

For mechanisms and kinetics, string methods approximate minimum (free-)energy paths / isocommittor-like structures in CV space \cite{E_Ren_VandenEijnden2002String,MaraglianoEtAl2006StringCV}.
Transition Path Sampling (TPS) samples whole reactive trajectories via Monte Carlo in path space, while Transition Interface Sampling (TIS) uses interfaces to compute rates more efficiently \cite{DellagoEtAl1998TPS,vanErpMoroniBolhuis2003TIS}.
Splitting/branching approaches such as Forward Flux Sampling (FFS), Weighted Ensemble (WE), and Adaptive Multilevel Splitting (AMS) estimate rare-event probabilities/rates by replicating promising trajectories (often guided by a score/CV)\cite{AllenFrenkelTenWolde2006FFS,HuberKim1996WE,CerouGuyader2007AMS}.
Milestoning discretizes progress along interfaces (“milestones”) and stitches short trajectory statistics to recover long-time kinetics \cite{FaradjianElber2004Milestoning}.

%%%%%%%%%%%%%%%%%%%%%%%%%%%%%%%
%%%%%%%%%%%%%%%%%%%%%%%%%%%%%%%
%%%%%%%%%%%%%%%%%%%%%%%%%%%%%%%
%%%%%%%%%%%%%%%%%%%%%%%%%%%%%%%
\section{Numerical Experiments}\label{sec:numexp}
%%%%%%%%%%%%%%%%%%%%%%%%%%%%%%%
%%%%%%%%%%%%%%%%%%%%%%%%%%%%%%%

In the subsequent, we will study the proposed lifting strategies for a simple 2d-system where clear interpretation and inspection of all desired properties is still intuitively possible.  In Sec.~\ref{sec:high-d test} of the Appendix, we show how to generalize this system into higher dimensions for more realistic numerical experiments.

\subsection{Simple test system}
We consider the 2-dimensional potential with two main wells, 
\begin{equation}\label{dw}
V_{dw}(x_1,x_2)=\alpha (x_1^2-1)^2+\beta (x_2^2-1)^2+ (1-\exp(-\gamma(x_1-x_2)^2))
\end{equation}
and consider the uncontrolled 2-dim diffusion in $x=(x_1,x_2)$ with $b=-\nabla V_{dw}$,
\begin{equation}\label{dyn-dw}
\frac{dx_t}{dt}=-\nabla V_{dw}(x) + \sigma \frac{dW_t}{dt}.
\end{equation}
In the following, we consider $\alpha=\beta=1$, $\gamma=2$, and $\sigma=0.7$.
Let $\mathcal{L}_{dw}$ denote the associated generator,
\begin{equation}\label{cL}
\mathcal{L}f(x)=-\nabla V_{dw}(x)\cdot\nabla f(x)+\frac{\sigma^2}{2}\Delta f(x),
\end{equation}
let $\varphi_{dw}$ be the eigenfunction of $\mathcal{L}_{dw}$ for its second eigenvalue $\lambda<0$ and assume that $\varphi_{dw}$ is bounded on the state space so that its minimum and maximum exist. 
Then, the membership function
\[
\chi(x_1,x_2)=\frac{\varphi(x_1,x_2)-\min\varphi}{\max\varphi-\min\varphi},
\]
is an appropriate, dynamically informative CV for the system, see \cite{SDWS25}. We will solely consider the CV $xi=\chi$ in this section.

%%%%%%%%%%%%%%%%%%%%%%
\subsection{Transition kinetics of the test system}
\begin{figure}
    \centering
    \includegraphics[width=\linewidth]{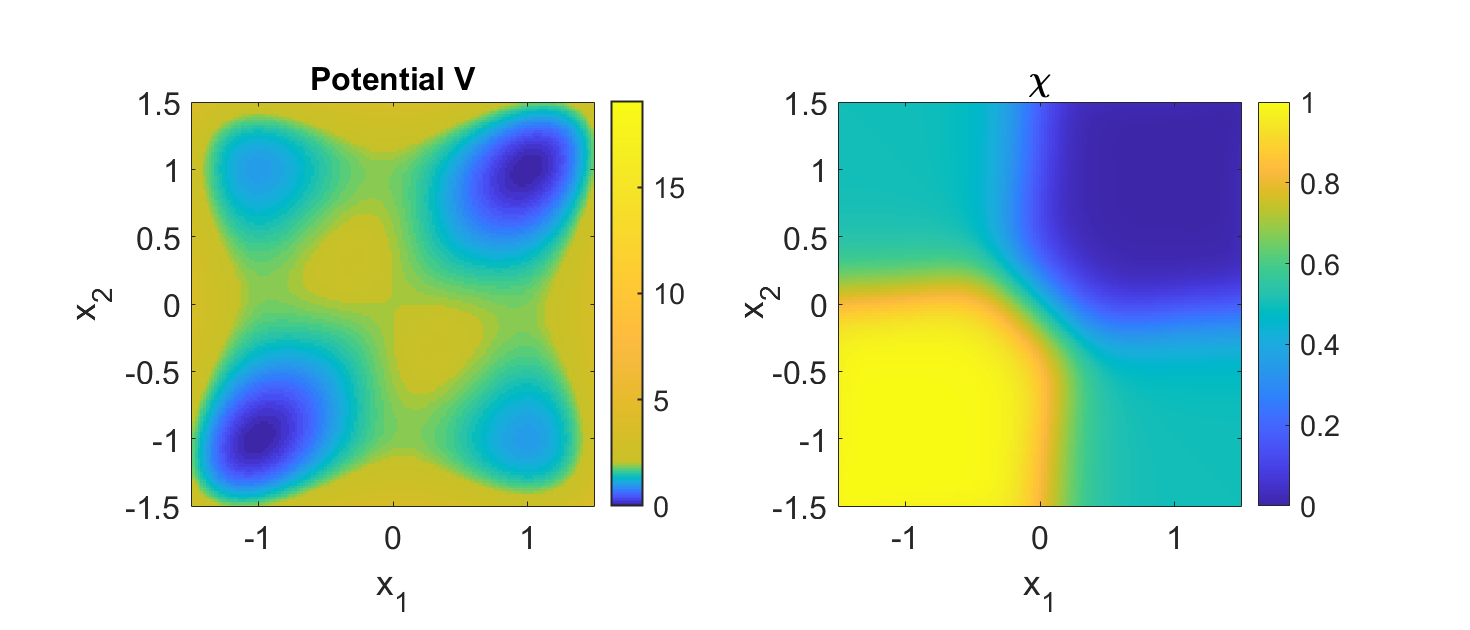}
    \vspace{-9mm}
    \caption{Illustration of potential function $V=V_{dw}$ with two main wells around $(x_1,x_2)=(-1,-1)$ as well as $(x_1,x_2)=(1,1)$, and associated $\chi$-function, associated with a dominant second eigenvalue $\lambda=-2.4\cdot 10^{-3}$ of $\mathcal{L}$ (with third eigenvalue $-6.6\cdot 10^{-3}$).}
    \label{fig:pot_chi}
\end{figure}

The generator of the 2d test system can be discretized via available techniques, based on finite elements or differences with no-flux boundary conditions in a finite-size discretization box, or via the Square-Root Approximation (SqRA) \cite{Lie2013, Donati2018b, Donati2021}. Based on these techniques, we CV $\xi=\chi$ can be computed, see Fig.~\ref{fig:pot_chi}.

\begin{figure}[h]
    \centering
    \includegraphics[width=0.35\linewidth]{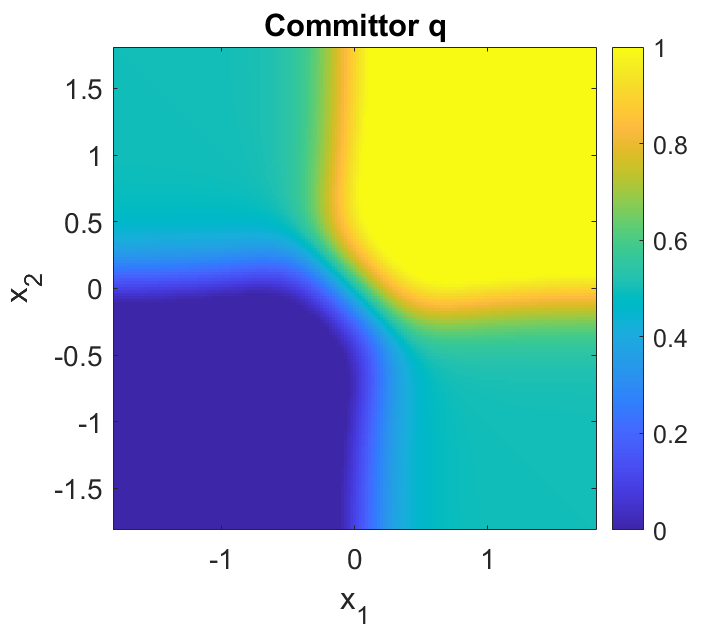}
    \includegraphics[width=0.34\linewidth]{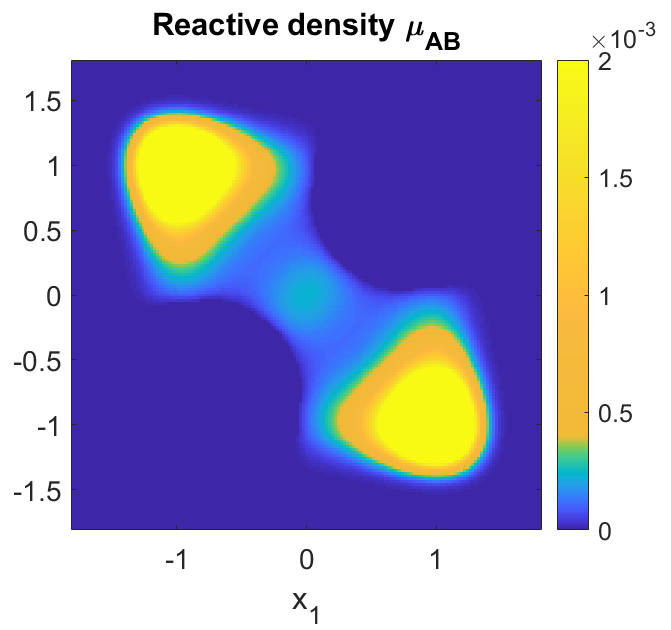}
    \includegraphics[width=0.29\linewidth]{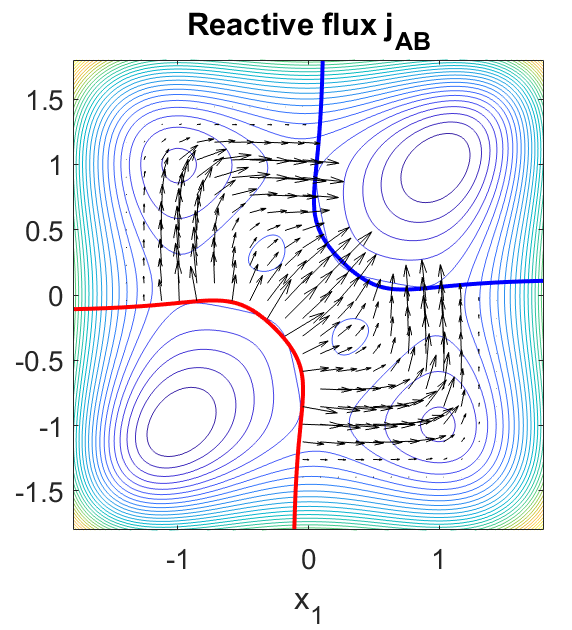}
    \vspace{-9mm}
    \caption{Full system committor function $q$ (left), reactive density $\mu_{AB}$ (middle), and reactive flux $j_{AB}$ for sets $A$ and $B$ as given in the text.}
    \label{fig:reactiveDensity}
\end{figure}

After a discretization of $\cL$ is available, it can be used to compute the committor function and further TPT quantities for the test system. Fig.~\ref{fig:reactiveDensity} shows the full system committor function $q$ for the sets $A=\{x: \chi(x)\le 0.1\}$ and $B=\{x: \chi(x)\ge 0.9\}$, the associated reactive density $\mu_{AB}$ (middle panel), and the associated reactive flux. 

\begin{figure}
    \centering
    \includegraphics[width=0.4\linewidth]{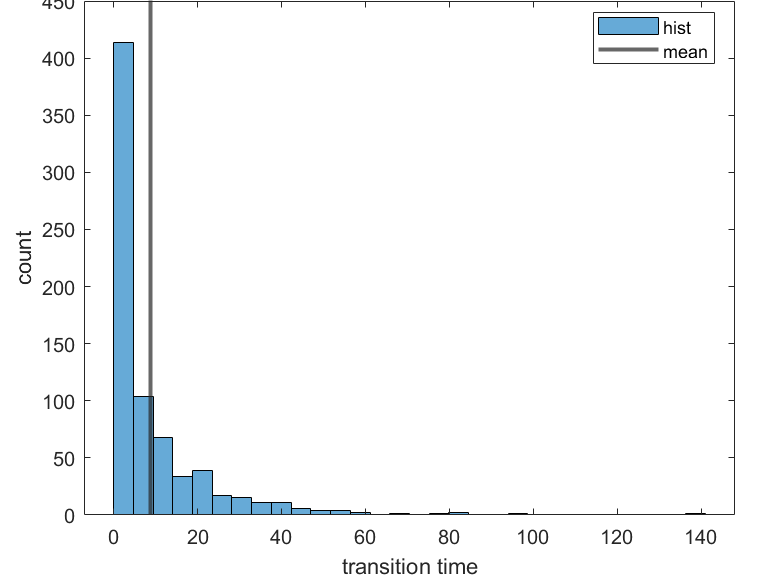}
    \caption{Statistics of transition times of the 2d test system for $\alpha=\beta=1$, $\gamma=2$ and $\sigma=0.7$. The statistics results from an ensemble of more than 600 reactive trajectories starting in $A=\{x: \chi(x)\le 0.1\}$ and going next to $B=\{x: \chi(x)\ge 0.9\}$ that were cut out a long-term simulation of the system.}
    \label{fig:transition_times_full}
\end{figure}

Another form of characterization of the transition kinetics of the test system is available by performing a very long simulation of the full system and cutting out all the reactive trajectories. Fig.~\ref{fig:transition_times_full} shows the resulting statistics for the lengths of these reactive trajectories, resulting in a average length of about $T=9.0$, which is an approximation of the expected transition length.

%%%%%%%%%%%%%%%%%%
\subsection{Effective dynamics}
According to \cite{SDWS25}, the effective dynamics for a system with CV $\xi=\chi$ has the form of a scalar Ornstein-Uhlenbeck process with $z$-dependent noise,
\begin{equation}\label{eff-dyn}
dz_t = (c+\lambda z_t)dt +\hat{\sigma}(z_t) dW_t,
\end{equation}
with constant $c=\lambda \min\varphi/(\max\varphi-\min\varphi)$, and 
\[
\hat{\sigma}(z)=\sigma\left(\mathbb{E}_{\hat{\mu}}[(\nabla\chi^\top \nabla \chi)(x)\;|\;\chi(x)=z]\right)^{1/2} =\sigma \left(\mathbb{E}_{\mu}[\|\nabla \chi\|^2(x_1,x_2)\;|\;\chi(x_1,x_2)=z]\right)^{1/2}
\]
where $\hat{\mu}(x) \propto\exp(-\beta V(x))$, with $\beta=2/\sigma^2$, $\|\cdot\|$ denotes the 2-norm in $\mathbb{R}^2$, and $\mu(x_1,x_2)\propto\exp(-\beta V_{dw}(x_1,x_2))$. According to \cite{Sikorskietal2025}, the effective dynamics is characterized by the resulting effective diffusion coefficient 
\[
D_{\text{eff}}(z)=\frac{1}{2}\hat{\sigma}(z)^2,
\]
and the effective potential
\[
V_{\text{eff}}(z)=\log(D_{\text{eff}}(z)) + \int_{z_\ast}^z \frac{c+\lambda z'}{D_{\text{eff}}(z')}dz',
\]
with arbitrary $z_\ast$ chosen such that $\min V_{\text{eff}}=0$.

\begin{figure}
    \centering
    \includegraphics[width=0.7\linewidth]{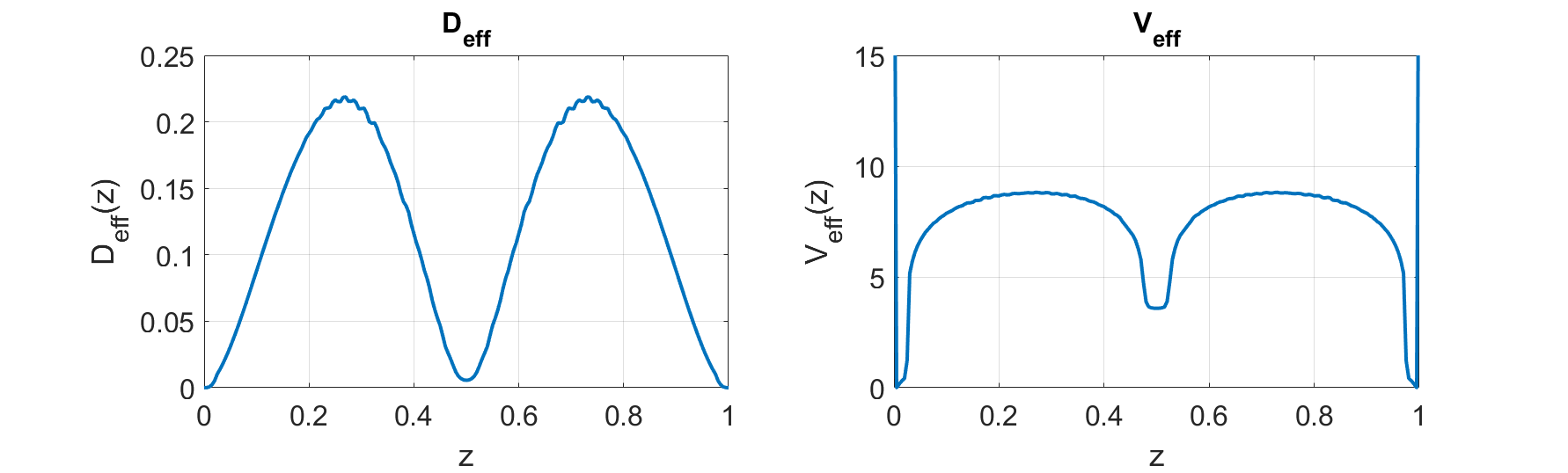}
    \vspace{-5mm}
    \caption{Effective diffusion coefficient $D_{\text{eff}}$ and effective potential $V_{\text{eff}}$ of the effective dynamics for the simple test system.}
    \label{fig:effpot}
\end{figure}

For the case of the simple test system, we get $c=0.0012$ and $\lambda=-0.0024$, and
Fig.~\ref{fig:effpot} illustrates the resulting effective diffusion coefficient 
and the effective potential.

\begin{figure}[h]
    \centering
    \includegraphics[width=0.8\linewidth]{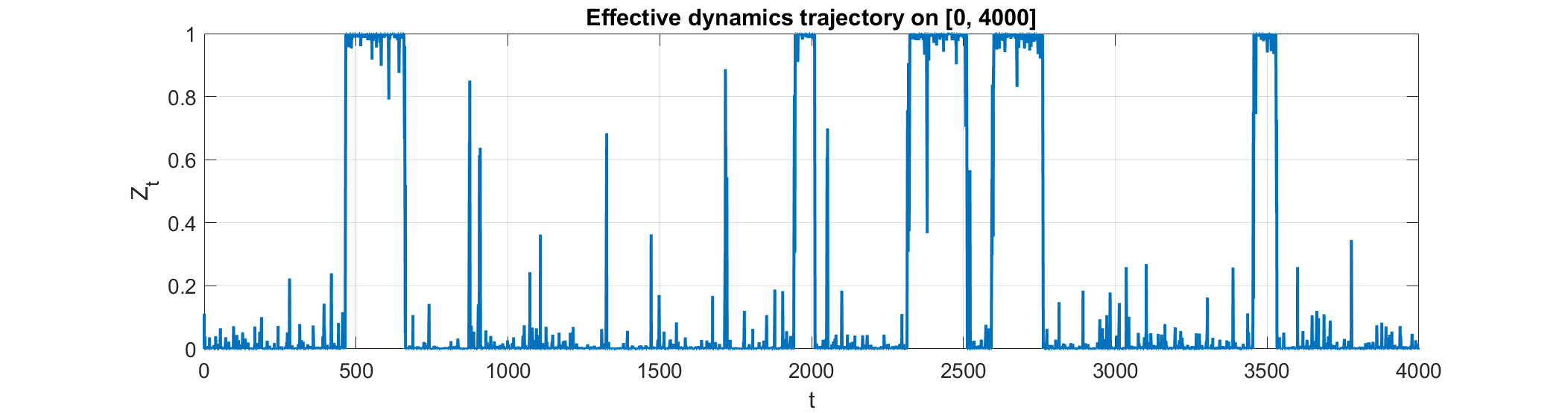}
    \vspace{-5mm}
    \caption{Typical trajectory of the effective dynamics, started in $z_0=0.05$.}
    \label{fig:efftraj}
\end{figure}

Integration of the SDE (\ref{eff-dyn}) of the effective dynamics in time yields trajectories that exhibit metastable dynamics in the 1-dim latent space $[0,1]$. As illustrated in Fig.~\ref{fig:efftraj} a typical trajectory fluctuates in the vicinity of the main metastable regions at $z=0$ and $z=1$ with rare transitions between them. We also observe a much weaker but visible metastable region around $z=0.5$. Computation of the associated effective generator $\mathcal{L}_{\text{eff}}$ is difficult because of the sensitivity wrt the behavior of $D_{eff}$ at the boundary $z=0$ and $z=1$. However, we can approximate  $\mathcal{L}_{\text{eff}}$ well by computing the associated Koopman operator $\mathcal{K}_\tau=\exp(\tau \mathcal{L}_{\text{eff}})$ via data from a long-term simulation of the effective dynamics, cf. \cite{SchuetteKlusHartmann2023}. Using $\tau=2$, a simulation of length $T=500.000$, and 200 uniform discretization boxes, we find 
 leading eigenvalues $\lambda_1=0$, and $\lambda_2=-0.0025$, showing that the dominant metastability of the full-dim system is well-reproduced. 

%%%%%%%%%%%%%%%%%%%%%%%%%%%%%%%%%%%%%%%%%%%%
%%%%%%%%%%%%%%%%%%%%%%%%%%%%%%%%%%%%%%%%%%%%
%%%%%%%%%%%%%%%%%%%%%%%%%%%%%%%%%%%%%%%%%%%%

\subsection{Guided diffusion bridges for the test system}
\begin{figure}[h]
    \centering
    \includegraphics[width=0.5\linewidth]{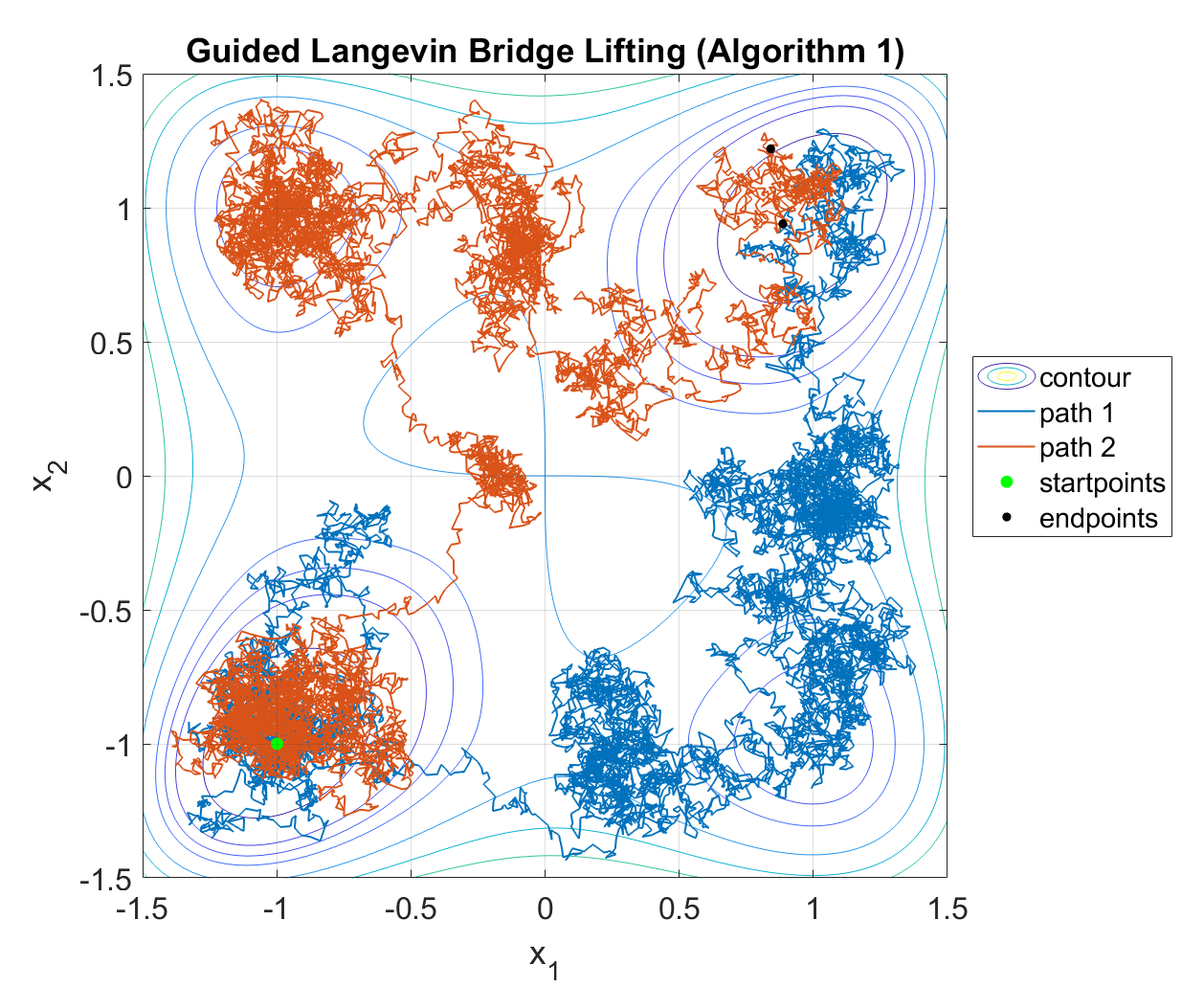}
    \vspace{-9mm}
    \caption{Results of Algorithm 1 for the 2-dim test system with potential $V_{dw}$ from (\ref{dw}), $\sigma=0.7$, and CV $\chi$ as displayed in  Fig.~\ref{fig:pot_chi} with 1-dim latent space $[0,1]$. Displayed are two typical paths resulting from Algorithm 1 on top of contour lines of the potential showing the main wells and the two side wells. }
    \label{fig:langevin_bridge}
\end{figure}

In a first test of Algorithm 1 we choose the information $Z_{T+j\Delta}$ given for the effective dynamics as a simple linear connection between $z\approx 0$ and $z\approx 1$. The results are shown in Fig.~\ref{fig:langevin_bridge}. The hyperparameters used for the algorithm were: $T=0$, $T_{new}=10$, gain schedule $G_t=100$, startpoint $X_0=(-1,-1)$, and $Z_{j\Delta t}=\chi(X_0)+j\delta z$ with $\delta z=(\chi(X_{10})-\chi(X_0))/10$, $\Delta t=1$, and trageted endpoint $X_{10}=(+1,+1)$. As expected we observe how the sampled diffusion bridge for the full dynamics, starting in $X_0=(-1,-1)$, slowing transits the transition region and ends up in the vicinity of $(+1,+1)$. 

%\vspace{-1cm}

\begin{figure}[h]
    \centering
    \includegraphics[width=0.5\linewidth]{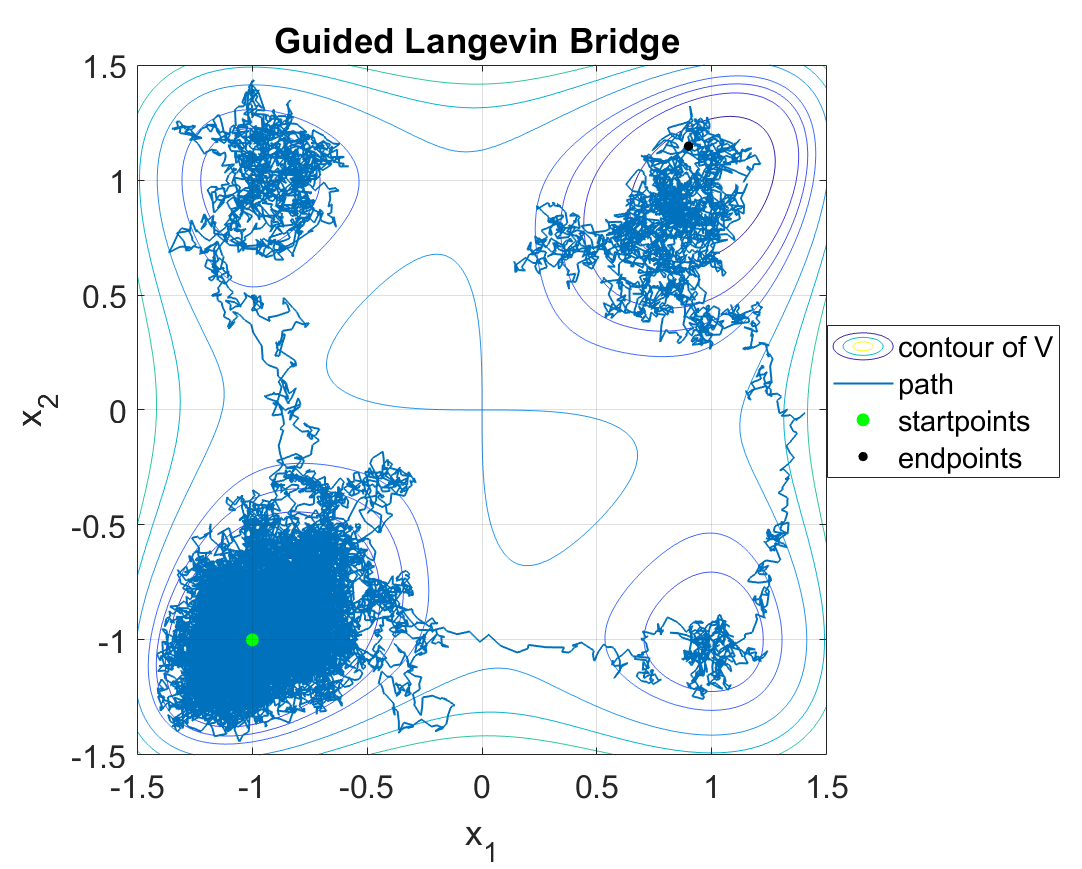}
    \vspace{-5mm}
    \caption{Results of Algorithm 1 for the 2-dim test system for the setting described in Fig.~\ref{fig:langevin_bridge}. Displayed is a typical paths driven by the effective dynamic on top of contour lines of the potential. }
    \label{fig:diffusion_bridge}
\end{figure}

Next, we repeat the numerical experiment keeping the entire setting but with $(Z_{j\Delta t})_{j=0,\ldots,400}$, $\Delta t=1$, resulting from a trajectory of the effective dynamics that is taken from the trajectory shown in Fig.~\ref{fig:efftraj}. More precisely, $(Z_{j\Delta t})$ is identical to the first part of the trajectory shown there in the initial time interval $[0,400]$. In Fig.~\ref{fig:diffusion_bridge}, we observe that the diffusion bridge starts in the main well around $(-1,-1)$ and fluctuates there first, then makes a transition to the side minimum around $(-1,+1)$, returns to the first main well and then proceeds towards the other main well around $(+1,+1)$ via the other side well around $(+1,-1)$.

%%%%%
\subsection{Transition path ensembles}

\begin{figure}[h]
    \centering
    \includegraphics[width=0.48\linewidth]{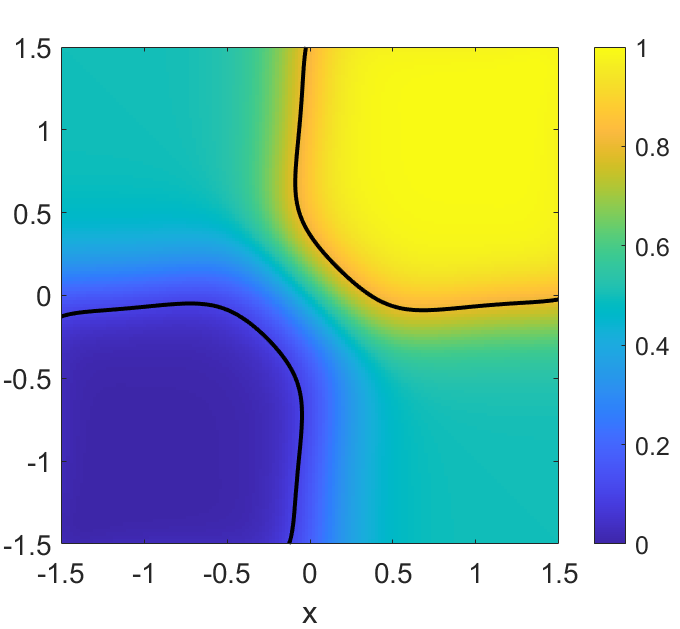}
    \includegraphics[width=0.46\linewidth]{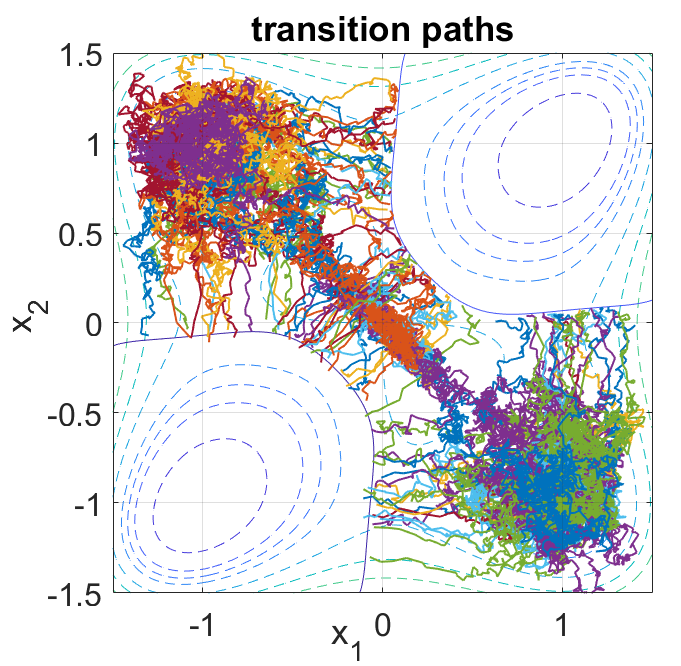}
    \vspace{-2mm}
    \caption{Results of Algorithm 1 for the 2-dim test system for the setting described in Fig.~\ref{fig:langevin_bridge}. Left: Illustration of the levelsets $\mathbb{L}_{0.1}$ and $\mathbb{L}_{0.9}$ on top of a colormap of the $\chi$-function. Right: Ensemble of reactive paths that start in $\mathbb{L}_{0.1}$ and hit $\mathbb{L}_{0.9}$ without returning to $\{z:\;z\le 0.1\}$, computed via bridge sampling based on a coarse transition path going from $z=0.1$ to $z=0.9$. }
    \label{fig:transition_path_ensemble}
\end{figure}

We can also use Algorithm 1 for sampling possible full-dimension \emph{reactive} trajectories, i.e., such trajectories that start on the levelset $\mathbb{L}_{z_{\text{min}}}$  and make the transition to the levelset of $z=z_{\text{max}}$, $0\le z_{\text{min}}<z_{\text{max}}<1$ without returning to $\mathbb{L}_{z_{\text{min}}}$ in between. Fig.~\ref{fig:transition_path_ensemble} shows such an ensemble of 100 reactive trajectories guided by a piece of the coarse trajectory displayed in Fig.~\ref{fig:efftraj} that directly goes from $z_{\text{min}}=0.1$ to $z_{\max}=0.9$. The trajectories shown in Fig.~\ref{fig:transition_path_ensemble} (right panel) are computed using a constant gain schedule $G_t=25$ that leads to strong guiding.  

\begin{figure}
    \centering
    \includegraphics[width=0.47\linewidth]{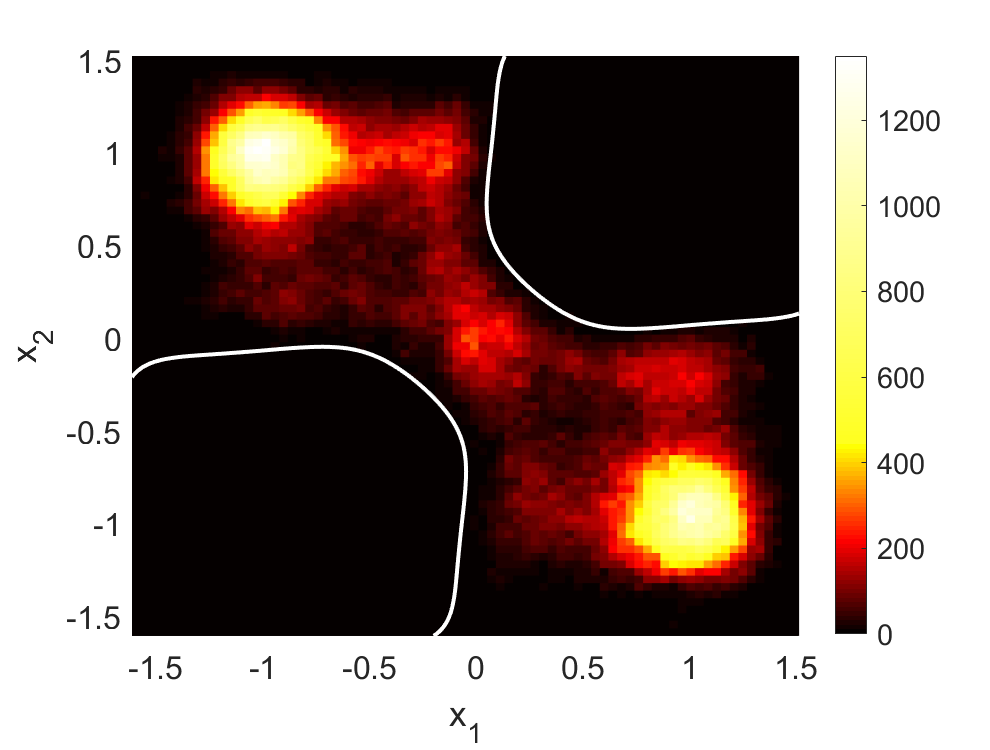}
    \includegraphics[width=0.45\linewidth]{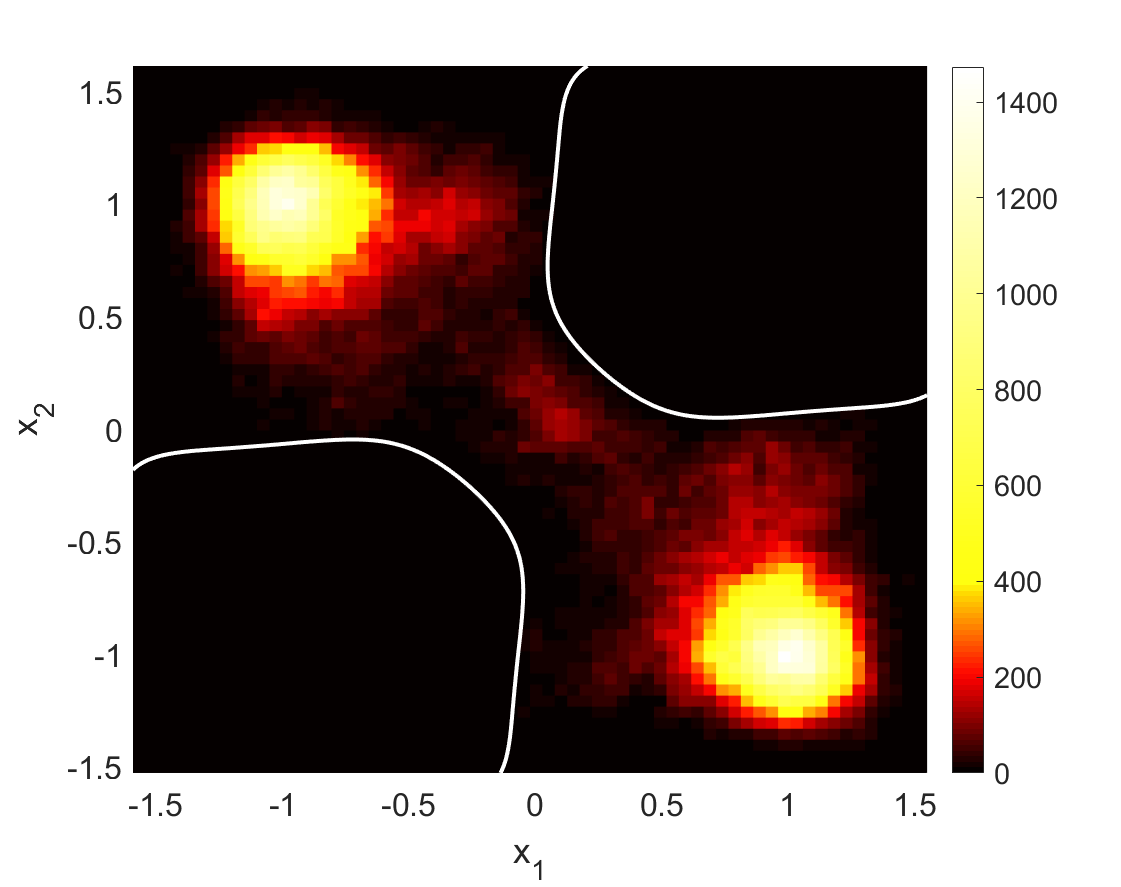}
    \caption{Histograms induced by 100 reactive trajectories computed with Algorithm 1. Left: Constant gain schedule $G_t=25$. Right: Small adaptive $G_t$. }
    \label{fig:histo_reactive}
\end{figure}

In Fig.~\ref{fig:histo_reactive} the effect of the gain schedule is illustrated: In the left panel the histogram induced by the ensemble of reactive trajectories of Fig~\ref{fig:transition_path_ensemble} with large and constant $G_t=25$  is shown. In the right panel, the respective histogram is plotted for the ensemble of reactive trajectories computed with small adaptive gain schedule $G_t$ is shown.   

These results must be compared to the correct full-system reactive density $\mu_{AB}$ given by TPT, see Fig.~\ref{fig:reactiveDensity} above. We observe good agreement with the reactive histogram computed by means of Algorithm 1. 

\begin{figure}
    \centering
    \includegraphics[width=0.3\linewidth]{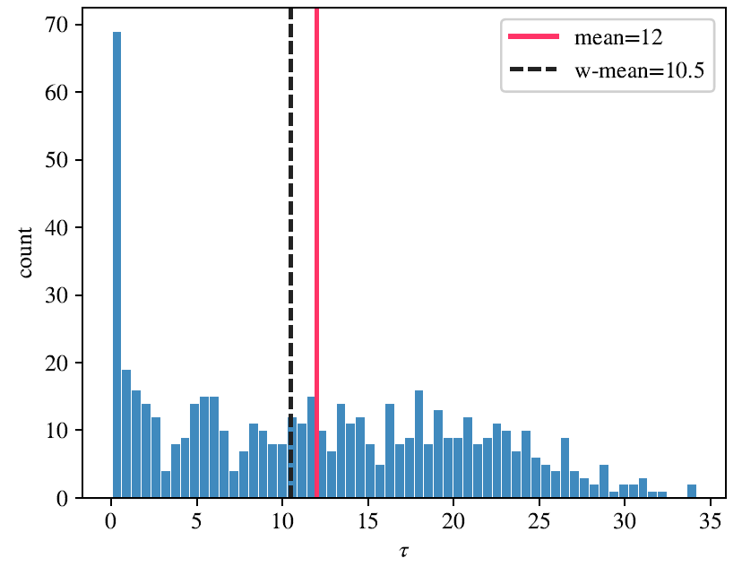}
    \includegraphics[width=0.3\linewidth]{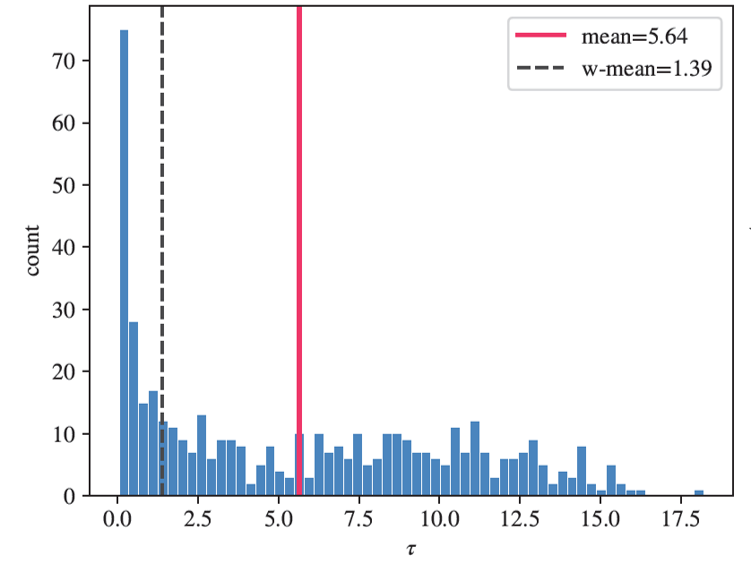}
    \includegraphics[width=0.3\linewidth]{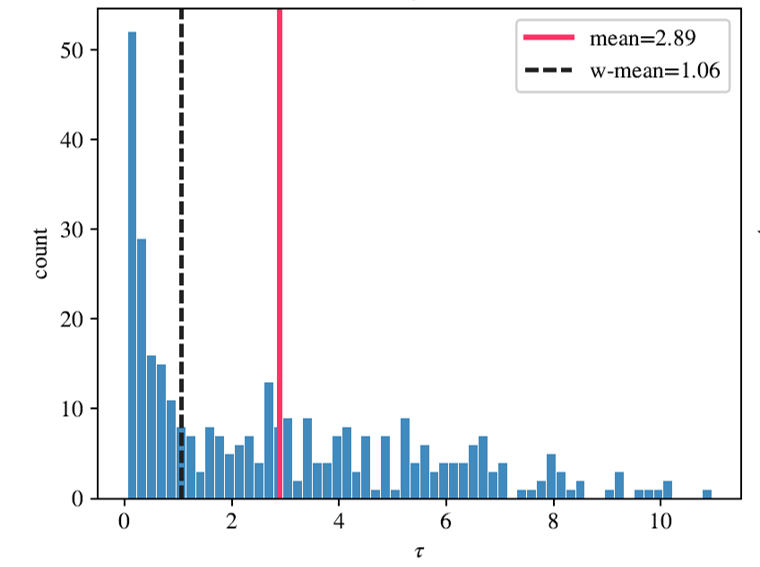}
\vspace{-5mm}    
\caption{statistics of the length of reactive pieces computed via Algorithm 1 for different choices of the gain schedule $G_t$. Left: $G_t=15$. Middle: $G_t=25$. Right: $G_t=50$. The computed average length changes from 12.0, via 5.6, to $2.9$.}
    \label{fig:stats_react_effdyn}
\end{figure}

However, Fig.~\ref{fig:stats_react_effdyn} shows that the statistics of the length of reactive pieces computed via Algorithm 1 are rather sensitive to the choice of the gain schedule $G_t$. 

\subsection{Transition probabilities}
Next, we turn to compute (small) transition probabilities of the form
\[
p_B(t,x)=\mathbb{P}(X_t\in B|X_0=x)
\]
for the uncontrolled diffusion $(X_t)$ of the full system for sets of the form $B=\{x:\; \chi(x)>z_\ast\}$. To this end, we computed 5.000 trajectories of the uncontrolled system, starting in $x_0=(-0.2,-0.2)$, and the Monte Carlo estimate $\hat{p}_B(T,x_0)$ for the probability $p(T,x_0)$ that $B$ with $z_\ast=0.9$ is reached at $T=20$, resulting in the MC estimate:
\[
p_B^{\text{MC}}(T,x_0) =0.148\pm 0.008.
\]
In order to prepare numerical experiments using the optimally guided diffusion we first computed the probability 
\[
p(s,z) =\mathbb{P}(Z_s>z_\ast|Z_0=z)
\]
for the (uncontrolled) effective dynamics $(Z_s)$ by solving the PDE (\ref{eq:bk-effective}) using the numerical approximation of the generator $\cL_{\text{eff}}$ as computed above. Fig.~\ref{fig:p(s,z)} shows the result, together with the resulting derivative $\partial_z \log p(s,z)$.

\begin{figure}
    \centering
    \includegraphics[width=\linewidth]{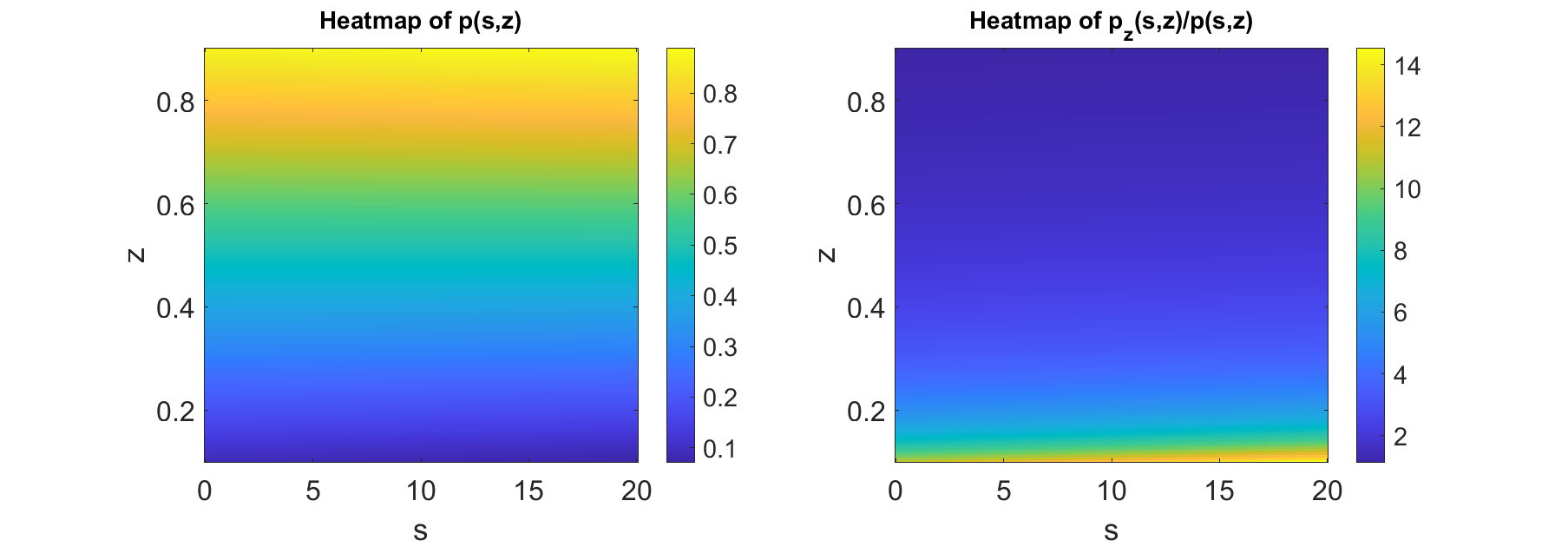}
    \vspace{-8mm}
    \caption{Heatmaps for probability $p(s,z)$ as defined in the text (left panel) and its derivative $\partial_z\log p(s,z)$ (right panel).}
    \label{fig:p(s,z)}
\end{figure}

Based on the computed $p(s,z)$ we started the optimal guidance algorithm of Sec.~\ref{ssec:optGuidance} with start point $x_0=(-0.2,-0.2)$. With $\kappa=1$, not all of 100 trajectories reached the target set. Using the boosting factor $\kappa=1.6$, the algorithm safely constructs trajectories that reach $B=\{x:\chi(x)>0.9\}$ in time $T=20$. Based on 100 trajectories the resulting estimate for the transition probability $p_B(T,x_0)$ is
\[
\hat{p}_B(T,x_0)=0.151\pm 0.012.
\]
Comparing this estimate with the MC-estimate $p_B^{\text{MC}}(T,x_0)$ computed via 5.000 trajectories of the uncontrolled full-dimensional system, we find sufficient agreement. The guided estimate $\hat{p}_B(T,x_0)$ does not show zero variance - which is an obvious consequence of the fact that the control $\tu$ used in the optimal guidance scheme is not the optimal control $u^\ast$ of (\ref{eq:u-star}) but an approximation of it. However, the estimate $\hat{p}_B(T,x_0)$ is based on 50 times less computational effort than the MC-estimate $p_B^{\text{MC}}(T,x_0)$ with comparable accuracy.

\subsection{Committor function from guided paths}

\begin{figure}
    \centering
    \includegraphics[width=0.5\linewidth]{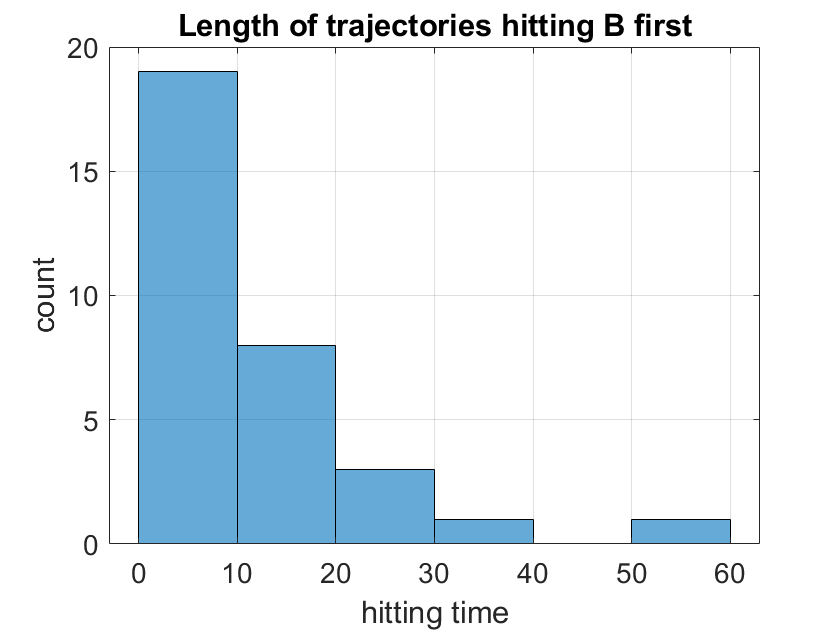}
    \caption{Histogram of hitting time $\tau_B$ of uncontrolled trajectories started in $x_0=(-1,0.2)$ that hit $B$ first.}
    \label{fig:length_committor}
\end{figure}

The committor function $q$ of the full system for the two sets $A=\{x:\,\chi(x)<0.1\}$ and $B=\{x:\,\chi(x)>0.9\}$ computed via the generator has already been shown in Fig.~\ref{fig:reactiveDensity} (left). For the point $x_0=(-1,0.2)$ we find $q(x_0)=0.3122$. When computing the committor based on uncontrolled trajectories of the full system via the direct Monte-Carlo estimate we find that, based on 100 such trajectories,
\[
\hat{q}^{\text{MC}}(x_0)= \frac{\text{no\, of\,trajectories\, hitting\,} B\,\text{first}}{\text{no\, of\,trajectories\, started}}=0.27 \pm 0.05.
\]
The average length of the trajectories that hit $B$ before $A$ is computed as $\tau_B=11.4\pm 12.1$; the histogram is shown in Fig.~\ref{fig:length_committor}.

\begin{figure}
    \centering
    \includegraphics[width=0.45\linewidth]{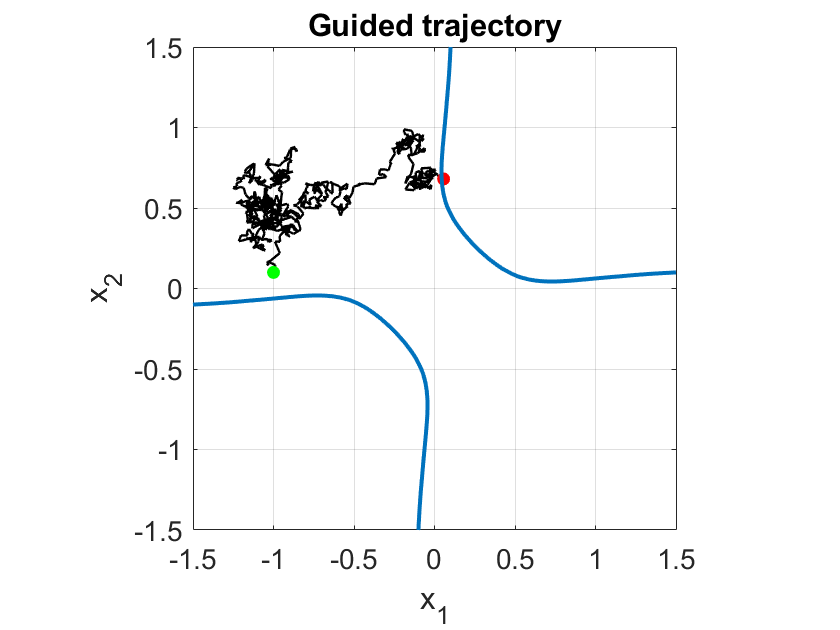}
    \includegraphics[width=0.45\linewidth]{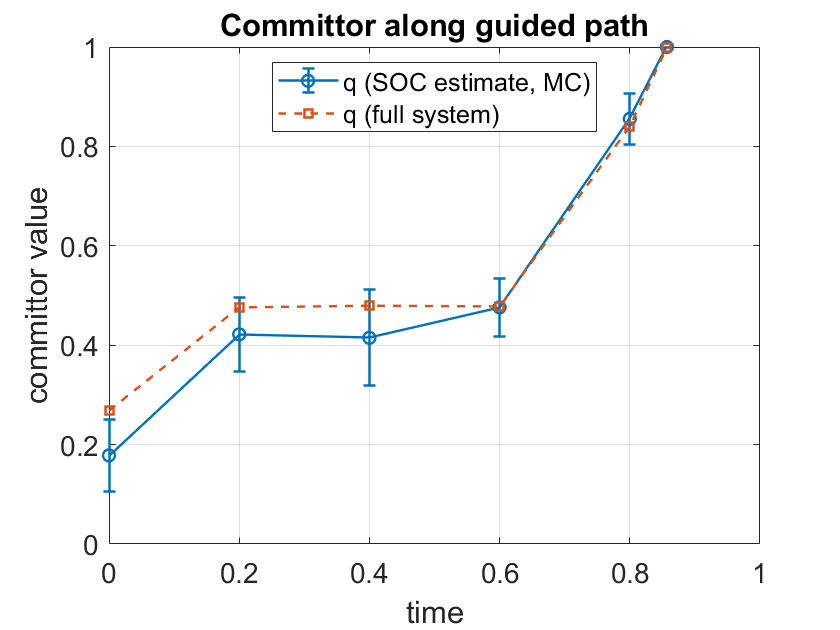}
    \caption{Typical guided trajectory for committor computation starting in $x_0=(-1,0.2)$ and hitting $B$ (left) and the committor $q$ and its guided estimate $\hat{q}$ along this trajectory. }
    \label{fig:committor_guided}
\end{figure}

When considering the estimate $\hat{q}$ given in (\ref{eq:committor_controlled}), based on the guidance scheme with  $\kappa=1.3$, started in $x_0=(-1,0.2)$ and using 100 guided trajectories, we find $\hat{q}(x_0)=0.26\pm 0.05$ with an average length $\tau_B=0.90\pm 0.11$ per guided trajectory, meaning that the guided estimate $\hat{q}$ is based on 10-fold less numerical effort than the uncontrolled MC estimate $\hat{q}^{\text{MC}}$ with comparable accuracy.  Fig.~\ref{fig:committor_guided} shows one of the guided trajectories (left) and the value of the committor $q$ and its estimate $\hat{q}$ along this trajectory.

%%%%%%%%%%%%%%%%%%%%%%%%%%%%%%%%%%%%%%%%%%%%%%%%%%%
%%%%%%%%%%%%%%%%%%%%%%%%%%%%%%%%%%%%%%%%%%%%%%%%%%%
\subsection{Alanine dipeptide in vacuum}

\begin{figure}
    \centering
    \includegraphics[width=0.35\linewidth]{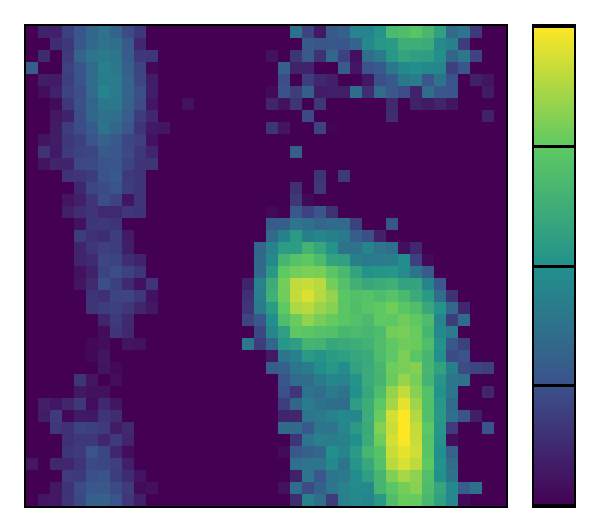}
    \includegraphics[width=0.35\linewidth]{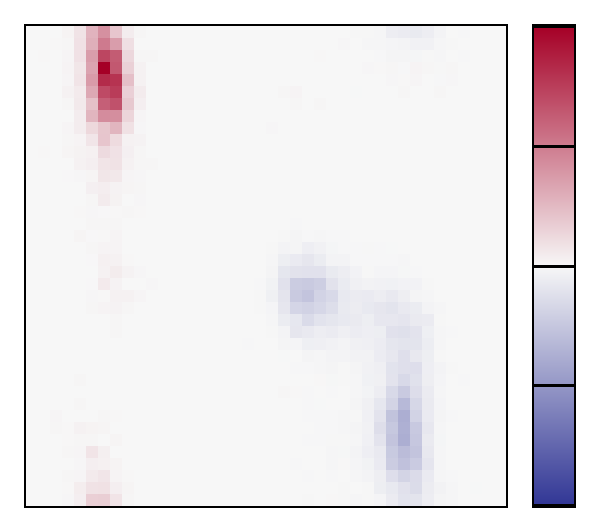}
    \caption{Illustration of the two dominant eigenmodes of the transfer operator $P$ (left for $\lambda=1.00$, right for $\lambda=0.91$ with corresponding colorbars from dark blue (almost 0 but positive) to yellow (largest positive), and from blue (negative) to red (positive). }
    \label{fig:P_eigv}
\end{figure}

As a first molecular example, we consider alanine dipeptide (ADP) in full atomsitic resolution in vacuum at $450\,\mathrm{K}$  and choose as collective variable the two-dimensional map
\[
\xi(x) = (\phi(x),\psi(x)),
\]
given by the two peptide dihedral angles. In contrast to the simple two-dimensional test problem,
the full dynamics now evolves in the atomistic configuration space, while the guidance is designed
on the low-dimensional Ramachandran domain. In the spirit of Sec.~4.3, we therefore combine a
coarse transition-probability computation in $(\phi,\psi)$-space with guided full-dimensional
simulation. More precisely, the Ramachandran domain is discretized by a periodic $40\times 40$
grid with boxes $A_i$. The corresponding $\tau=5\,\mathrm{ps}$ transfer operator $P$ is a $40^2\times 40^2$ square matrix whose entries are approximate transition probabilities $p(\tau,A_i,A_j)$ computed from 100 trajectories started in each box. Its dominant eigenvalues are $1.00$ and $0.91$ followed by $0.11$, where the eigenmodes of the first two are displayed in Fig.~\ref{fig:P_eigv} who clearly separate the two main conformations of ADP in the lower right (main conformation) and the upper left (side conformation) parts of the Ramachandran domain. 

\begin{figure}[htbp]
  \centering
  \includegraphics[width=0.4\textwidth]{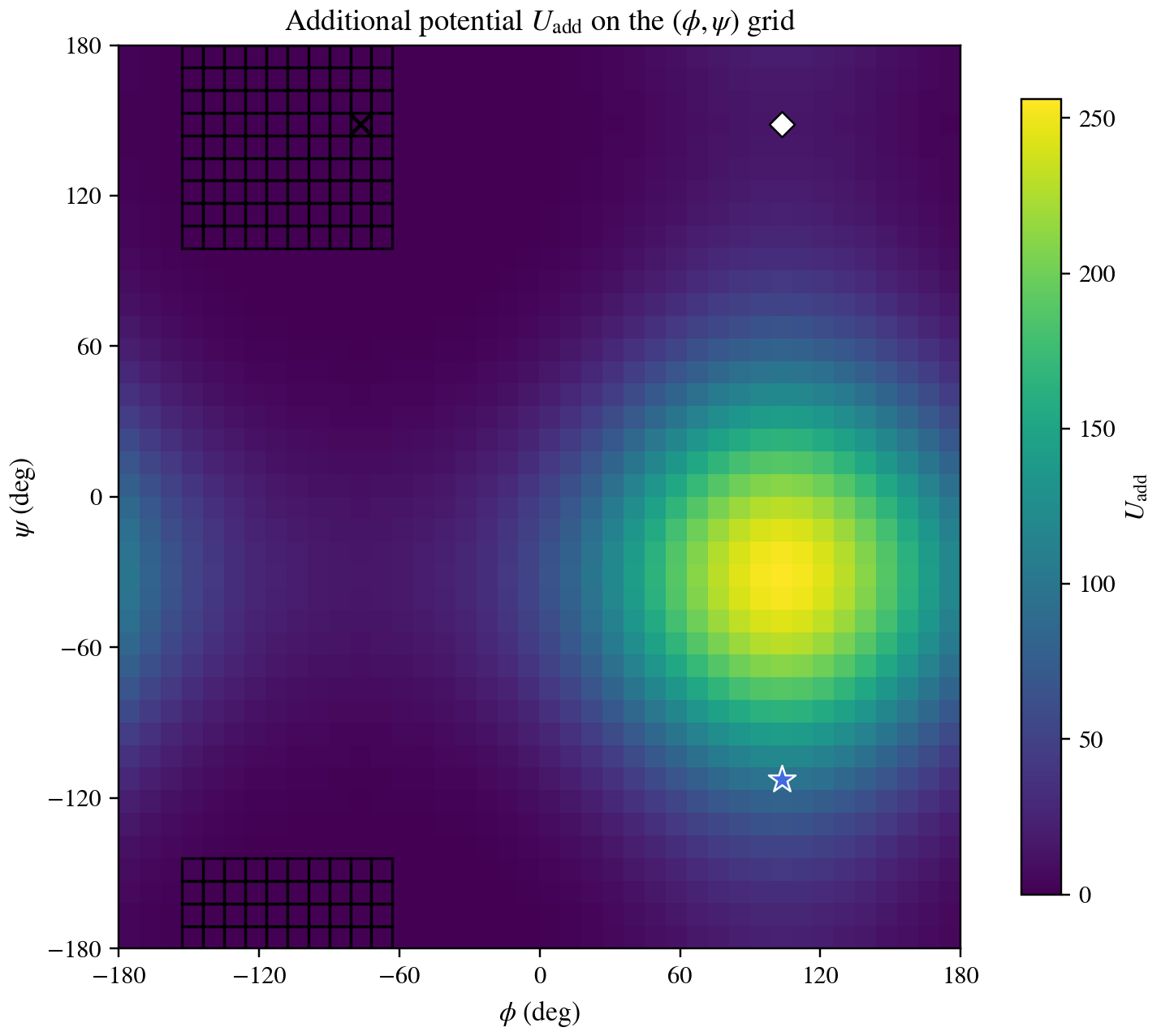}
  \caption{The quartic periodic added potential on the Ramachandran grid, together with the target set $B$ (black cluster of boxes).}
  \label{fig:uadd}
\end{figure}

Based on $P$ compute the boxwise finite-time transition probability
\[
p(s_k,A_i,B)=\bigl(P^{N-k}\mathbf{1}_B\bigr)_i,\quad N=T/\Delta t=20,\quad \Delta t=5\,\mathrm{ps},\quad s_k=k\Delta t,
\]
for the horizon $T=100\,\mathrm{ps}$, where the target set $B$ is chosen as a rectangular region in the upper-left sector of the
  $(\phi,\psi)$-plane, centered at $(-108.0^\circ,157.5^\circ)$, see Fig.~\ref{fig:uadd}. As in Sec.~4.3, this yields the probability-based part of the guidance through the discrete gradient of $\log p(s_k,A_i,B)$. Since in the present molecular
example $\xi$ is two-dimensional, the guidance acts in the $(\phi,\psi)$ variables and is then
realized in the full atomistic dynamics through the Jacobian of $\xi$, exactly in the form discussed in Secs.~3 and~4.3. 
Since the distribution $p(s,\cdot,B)$ is rather flat in the domain $\phi>0$ of the Ramachandran plane, the guiding force includes an additive part resulting from a smooth periodic auxiliary bias potential
  \[
  U_{\mathrm{add}}(\phi,\psi)
  =
  \bigl(1-\cos(\psi-\psi_0)+1+\cos(\phi-\phi_0)\bigr)^4,
  \qquad
  (\phi_0,\psi_0)=(103.5^\circ,148.5^\circ),
  \]
  whose minimum lies inside the target sector, see Fig.~\ref{fig:uadd}. Hence the guidance combines the transition-probability
  information from Sec.~4.3 with an additional forcing toward $B$ (via Girsanov weights).

In practice, OpenMM \cite{Eastman2023OpenMM8} is used as the engine for the full-dimensional molecular dynamics simulations, whereas the guidance field is computed on the Ramachandran grid and added to the physical dynamics as an additional biasing force acting through the CV map $\xi=(\phi,\psi)$. Thus the coarse model supplies the transition information, while the microscopic trajectory generation remains fully atomistic. 

 The transition probability $p_B(x_0,T)$ can roughly be estimated from the 5 ps transfer operator $P$ via $p(T,A,B)=(P^N\mathbf{1}_B)_\ell$, with $A=A_\ell$ denoting the unique set in which $x_0$ is located, which results in
  \[
  p_B(T,x_0)\approx 4.34\times 10^{-3}.
  \]
By simply starting 500 \emph{unbiased, unguided} trajectories in $x_0$ and counting the ones that end in the target $B$ at time $T$ = 500 ps, we have the direct Monte Carlo estimate 
\[
 p_B(T,x_0)\approx 1.8 \times 10^{-2},
\]
where the corresponding Wilson $95 \%$ interval is $[9.5 \times 10^{-3}, 3.4 \times 10^{-2}]$.

  We then generated an ensemble of $500$ \emph{guided} OpenMM trajectories. As in Sec.~4.3, the relevant
  quantity is a final-time occupancy probability rather than a first-hitting probability; nevertheless
  it is informative to record that, among the $500$ guided trajectories, $487$ visit $B$  before time $T$, and $463$ are in $B$ at the final time. Thus the guidance makes the transition to
  the target region typical under the controlled dynamics, while a small fraction of trajectories
  still leaves $B$ again before the terminal time.

 For this setup, the finite-time probability $p(x_0,T)$ from the prescribed initial condition $x_0$ is estimated by 
  \[
  \widehat p_{B,\varepsilon}(T,x_0)\approx 1.77\times 10^{-3}.
  \]
  A simple Monte Carlo error analysis gives a $95\%$ interval of order
  $[3\times 10^{-4},\,9\times 10^{-3}]$, so the estimate is clearly positive but still carries visible
  sampling uncertainty. Overall, this experiment shows that the Sec.~4.3 guidance principle can be
  transferred to a realistic molecular-dynamics setting with a two-dimensional CV and that,
  complemented by a smooth periodic bias, it yields a substantial enrichment of ADP transition paths
  into the prescribed target sector of the Ramachandran plane.

  Next, we repeat the same numerical experiment with discretized transfer operator with $\tau=2$ ps and $T=20$ ps, based on the identical box decomposition, target set and guidance set-up. We find that the second eigenvalue of the transfer operator now increases to $\lambda=0.96$ while the two leading eigenmodes still show the pattern illustrated in Fig.~\ref{fig:P_eigv}. The transition probability from $A$ to $B$ computed from the 2 ps transfer operator can again be taken as a rough estimator for the transition probability from $x_0$ to $B$, yielding
  \[
  p_B(T,x_0)\approx 2.6\times 10^{-5},
  \] 
  while the  estimate based on 500 unbiased, unguided trajectories yields $p_B(T,x_0)=0$ (none of the trajectories hit $B$).
 In comparison, the estimate based on 500 guided trajectories yields 
  \[
  \widehat p_{B,\varepsilon}(T,x_0)\approx 5.2\times 10^{-6},
  \]
  with a Monte Carlo error analysis based $95\%$ interval of order 
  $[9.3\times 10^{-7},\,1.0\times 10^{-4}]$.
  
%%%%%%%%%%
%%%%%%%%%%%%%%%%%%%%%%%%%%%%%%%%%%%%%%%%%%%%
%%%%%%%%%%%%%%%%%%%%%%%%%%%%%%%%%%%%%%%%%%%%
%%%%%%%%%%%%%%%%%%%%%%%%%%%%%%%%%%%%%%%%%%%%
%\clearpage
\section{Conclusion and Outlook}\label{sec:conclusion}

We studied the  lifting bottleneck  that arises when an informative collective variable (CV) and a reliable effective dynamics in CV space are available: coarse CV trajectories can be generated cheaply, but turning them into   dynamically consistent full-dimensional states and transition-path ensembles   is nontrivial without global sampling of conditional fiber measures. To address this, we proposed a   local, on-the-fly lifting strategy   that uses coarse information only as a  reference  and reconstructs microscopic realizations by simulating   guided full-system trajectories   that track the reference in CV space. 

Our first approach is tailored towards \emph{exploration} and proceeds by (i) interpolating a coarse CV trajectory from the effective dynamics to obtain a continuous reference path, (ii) generating an ensemble of guided full-dimensional trajectories by adding a feedback control based on the Jacobian of the CV map, and (iii) correcting the bias introduced by guidance through pathwise Girsanov reweighting. This guided bridge lifting yields a weighted ensemble of full trajectories and endpoints, and it retains the key locality property: it requires only evaluations of the drift, the CV, and its Jacobian along simulated paths. 

Beyond the exploration perspective, we connected the construction of guided trajectories to  stochastic optimal control (SOC) and the change-of-measure view of variance reduction. This clarifies how, for specific objectives (e.g., finite-time transition probabilities or committors), an optimal feedback control (for which accurate computation is infeasible) can be approximated using low-dimensional effective-dynamics information. This leads to guidance laws that make rare events more typical under the guidance and thereby reduce estimator variance. In this sense, the method provides both (a) a practical mechanism generator for full-dimensional transition paths driven by coarse predictions and (b) a principled entry point for designing variance-reducing guidance strategies informed by effective models. 

The numerical experiments on a metastable two-dimensional test system illustrate these roles clearly. Guided bridge lifting converts inexpensive coarse CV transition segments into   physically plausible full-system transition pathways, including barrier crossings and alternative channels (via side wells) when such behavior is present in the coarse trajectory. Ensembles of guided reactive paths reproduce key transition-path statistics (e.g., reactive density structure). Our experiments  also highlight a central practical trade-off: stronger gains improve CV tracking but can increase weight degeneracy, motivating diagnostics such as the effective sample size (ESS) and adaptive gain schedules. For quantitative tasks, the SOC-informed guidance demonstrates substantial efficiency gains: transition probabilities and committor-related quantities can be estimated with markedly reduced computational effort compared with direct Monte Carlo, while maintaining comparable accuracy in the reported tests. 

Several limitations and extensions follow naturally. The approach depends critically on CV quality: if the CV $\xi$ does not resolve the slow mechanism, aggressive guidance may yield unrealistic microscopic behavior and/or severe weight collapse. Robust long-horizon lifting suggests incorporating sequential Monte Carlo   (resampling/mutation over sub-intervals) and improved preconditioning or clipping strategies when the Jacobian of $\xi$ is ill-conditioned. Methodologically, a major direction is to move guidance closer to the optimal control by enriching the feedback with additional coarse information (e.g., policy iteration, local fiber geometry, or learned corrections). Future research will also have to work out extending the framework to applicability to more realistic molecular simulation settings. 

Overall, the lifting viewpoint developed here offers a practical route to   couple cheap coarse transition predictions with full-dimensional path generation. In combination with modern CV learning and effective-dynamics modeling, this provides a flexible building block for multiscale exploration and computation of rare transitions, capturing mechanistic variability through weighted ensembles while keeping the full-dimensional simulation effort focused where it matters most. 

\section*{Acknowledgment}
This research has been funded by Deutsche Forschungsgemeinschaft
(DFG) through grant \textit{CRC 1114 Scaling Cascades in Complex Systems} (Project
No. 235221301), Projects A05 and B03, and under Germany's Excellence Strategy MATH+: Berlin Mathematics Research Center (EXC 2046/1, project 390685689).

\bibliographystyle{unsrt}
\bibliography{refs}

@article{Eastman2023OpenMM8,
  author       = {Peter Eastman and Raimondas Galvelis and Ra{\'u}l P. Pel{\'a}ez and Charlles R. A. Abreu and Stephen E. Farr and Emilio Gallicchio and Anton Gorenko and Michael M. Henry and Frank Hu and Jing Huang and Andreas Kr{\"a}mer and Julien Michel and Joshua A. Mitchell and Vijay S. Pande and Jo{\~a}o PGLM Rodrigues and Jaime Rodriguez-Guerra and Andrew C. Simmonett and Sukrit Singh and Jason Swails and Philip Turner and Yuanqing Wang and Ivy Zhang and John D. Chodera and Gianni De Fabritiis and Thomas E. Markland},
  title        = {{OpenMM 8}: Molecular Dynamics Simulation with Machine Learning Potentials},
  journal      = {Journal of Physical Chemistry B},
  year         = {2023},
  volume       = {128},
  number       = {1},
  pages        = {109--116},
  doi          = {10.1021/acs.jpcb.3c06662}
}

@article{SDWS25,
author={Alexander Sikorski and Luca Donati and Marcus Weber and Christof Schuette},
year={2025},
title={Effective Dynamics and Transition Pathways from {K}oopman-Inspired Neural Learning of Collective Variables}, 
journal={Submitted to {Communications in Applied Mathematics and Computational Science}}
}

@article{eff-dyn-isolines,
author={Fatima-Zahrae Akhyar and Wei Zhang and Gabriel Stoltz and Christof Schuette},
year={2025},
title={Generative modeling of conditional probability distributions on the level-sets of collective variables}, 
journal={to be submitted to Communications in Applied Mathematics and Computational Science}
}

@article{HartmannSchuetteZhang2016Nonlinearity,
  author       = {Hartmann, Carsten and Sch{\"u}tte, Christof and Zhang, Wei},
  title        = {Model reduction algorithms for optimal control and importance sampling of diffusions},
  journal      = {Nonlinearity},
  year         = {2016},
  volume       = {29},
  number       = {8},
  pages        = {2298--2326},
  doi          = {10.1088/0951-7715/29/8/2298}
}

@article{Lie2013,
author = {Lie, Han Cheng and Fackeldey, Konstantin and Weber, Marcus},
journal = {SIAM. J. Matrix Anal. Appl.},
pages = {738–756},
title = {A Square Root Approximation of Transition Rates for a {M}arkov {S}tate {M}odel},
volume = {34},
year = {2013}
}

@article{Donati2018b,
  author={Donati, L. and Heida, M. and Keller, B. G. and Weber, M.},
    title={Estimation of the infinitesimal generator by square-root approximation},
  journal={J. Phys. Condens. Matter},
  year={2018},
  volume = {30},
  pages = {425201}
}

@article{Donati2021,
  author  = {Donati, L. and Weber, M. and Keller, B.~G.}, 
  title   = {Markov models from the square root approximation of the {F}okker–{P}lanck equation: {C}alculating the grid-dependent flux},
  journal = {J. Phys. Condens. Matter},
  year    = 2021,
  pages   = {115902},
  volume  = 33,
}

@article{SchuetteKlusHartmann2023,
  author    = {Christof Schütte and Stefan Klus and Carsten Hartmann},
  title     = {Overcoming the Timescale Barrier in Molecular Dynamics: Transfer Operators, Variational Principles, and Machine Learning},
  journal   = {Acta Numerica},
  volume    = {32},
  pages     = {517--673},
  year      = {2023},
  publisher = {Cambridge University Press},
  doi       = {10.1017/S0962492923000062},
}

@article{ZhangSchuette2025,
  author    = {Wei Zhang and Christof Schütte},
  title     = {On Finding Optimal Collective Variables for Complex Systems by Minimizing the Deviation Between Effective and Full Dynamics},
  journal   = {Multiscale Modeling \& Simulation},
  volume    = {23},
  number    = {2},
  pages     = {924--958},
  year      = {2025},
  doi       = {10.1137/24M1658917}
}

@article{ZhangHartmannSchuette2016,
  author    = {Wei Zhang and Carsten Hartmann and Christof Schütte},
  title     = {Effective dynamics along given reaction coordinates and reaction rate theory},
  journal   = {Faraday Discussions},
  volume    = {195},
  pages     = {365--394},
  year      = {2016},
  doi       = {10.1039/C6FD00147E}
}

@article{Hartmannetal2026,
title={Importance sampling with unbounded random
stopping times: computing committor functions
and exit rates without reweighting},
author={Carsten Hartmann and Annika Joester and Christof Schuette and  Alexander Sikorski and Marcus Weber},
journal={arXiv:2601.01489},
year={2026}
}

@article{HartmannSchuette2012OptimalNoneqForcing,
  author  = {Hartmann, Carsten and Sch{\"u}tte, Christof},
  title   = {Efficient rare event simulation by optimal nonequilibrium forcing},
  journal = {Journal of Statistical Mechanics: Theory and Experiment},
  year    = {2012},
  number  = {11},
  pages   = {P11004},
  doi     = {10.1088/1742-5468/2012/11/P11004}
}

@article{HartmannKebirieial2019,
    author = {Hartmann, Carsten and Kebiri, Omar and Neureither, Lara and Richter, Lorenz},
    title = {Variational approach to rare event simulation using least-squares regression},
    journal = {Chaos: An Interdisciplinary Journal of Nonlinear Science},
    volume = {29},
    number = {6},
    pages = {063107},
    year = {2019},
    month = {06},
    issn = {1054-1500},
    doi = {10.1063/1.5090271},
    url = {https://doi.org/10.1063/1.5090271},
    eprint = {https://pubs.aip.org/aip/cha/article-pdf/doi/10.1063/1.5090271/13534861/063107_1_online.pdf},
}

@book{SchuetteSarich2014,
author={Sch{\"u}tte, C. and Sarich, M.},
title={Metastability and {M}arkov State Models in Molecular Dynamics: {M}odeling, Analysis, Algorithmic Approaches}, 
series={Courant Lecture Notes No. 32}, 
publisher={American Mathematical Society},
year={2014}
}

@ARTICLE{BolhuisChandlerDellagoGeissler2002TPS,
  author = {P.~G. Bolhuis and C.~Dellago and D.~Chandler and P.~Geissler},
  title = {Transition path sampling: Throwing ropes over mountain passes, in
	the dark},
  journal = {Ann. Rev. of Phys. Chem.},
  year = {2001},
  note = {in press}
}

@ARTICLE{EVE2006JSP,
  author = {W. {E} and E. {Vanden-Eijnden}},
  title = {Towards a theory of transition paths},
  journal = {Journal of statistical physics},
  year = {2006},
  volume = {123},
  pages = {503-523}
}

@ARTICLE{MetznerSchuetteVandenEijnden2006TPTexamples,
  author = {P. Metzner and Sch{\"u}tte, C. and Vanden-Eijnden, E.},
  title = {Illustration of Transition Path Theory on a Collection of Simple
	Examples},
  journal = {J. Chem. Phys.},
  year = {2006},
  volume = {125},
  number = {8},
  note = {084110}
}

@article{EVE2010TPTreview,
author={W.~{E} and E.~{Vanden-Eijnden}},
title={Transition-path theory and path-finding algorithms for the study of rare events},
journal={Annu. Rev. Phys. Chem.}, 
volume={61},
pages={391--420},
year= {2010}
}

@article{BerezhkovskiiSzabo2019CommittorsMilestones,
author = {Berezhkovskii, A. M.  and Szabo, A.},
title = {Committors, first-passage times, fluxes, Markov states, milestones, and all that},
journal = {The Journal of Chemical Physics},
volume = {150},
number = {5},
pages = {054106},
year = {2019}
}

@ARTICLE{MoroniVanErpBolhuis2004TIS,
  author = {D. Moroni and T. {van Erp} and P. Bolhuis},
  title = {Investigating rare events by transition interface sampling},
  journal = {Physica {A}},
  year = {2004},
  volume = {340},
  pages = {395-401}
}

@ARTICLE{FaradjianElber2004Milestoning,
  author = {A. K. Faradjian and R. Elber},
  title = {Computing time scales from reaction coordinates by milestoning},
  journal = {J. Chem. Phys.},
  year = {2004},
  volume = {120},
  pages = {10880-10889}
}

@article{KlusNueskeKoltaiWuKevrekidisSchuetteNoe2018,
    author  = {S. Klus and F. N\"uske and P. Koltai and H. Wu and I. Kevrekidis and C. Sch\"utte and F. No\'e},
    title   = {Data-driven model reduction and transfer operator approximation},
    journal = {Journal of Nonlinear Science},
    volume  = {28},
    issue   = {3},
    pages   = {985--1010},
    year    = {2018},
    doi     = {10.1007/s00332-017-9437-7},
}

@book{BowmanPandeNoe2013,
  title={An introduction to Markov state models and their application to long timescale molecular simulation},
  author={Bowman, G. R and Pande, V. S and No{\'e}, F.},
  volume={797},
  year={2013},
  publisher={Springer Science \& Business Media}
}

@article{RabbenRayWeber2020ISOKANN,
author = {Rabben, R. J. and Ray, S. and Weber, M.},
title = {ISOKANN: {I}nvariant subspaces of {K}oopman operators learned by a neural network},
journal = {The Journal of Chemical Physics},
volume = {153},
number = {11},
pages = {114109},
year = {2020},
doi = {10.1063/5.0015132}
}

@article{Sikorskietal2025,
    author = {Alexander Sikorski and Luca Donati and Marcus Weber and Christof Schuette},
    title = {Effective Dynamics and Transition Pathways from {K}oopman-Inspired Neural Learning of Collective Variables},
    journal = {submitted to to Communications in Applied Mathematics and Computational Science},
    year = {2025}
}

@article{NueskeKoltaiBoninsegnaClementi2021,
  author    = {Feliks N{\"u}ske and P{\'e}ter Koltai and Lorenzo Boninsegna and Cecilia Clementi},
  title     = {Spectral Properties of Effective Dynamics from Conditional Expectations},
  journal   = {Entropy},
  volume    = {23},
  number    = {2},
  pages     = {134},
  year      = {2021},
  doi       = {10.3390/e23020134}
}

@article{LegollLelievreOlla2017,
  author    = {F. Legoll and T. Leli{\`e}vre and S. Olla},
  title     = {Pathwise estimates for an effective dynamics},
  journal   = {Stochastic Processes and their Applications},
  volume    = {127},
  year      = {2017},
  pages     = {2841--2863},
  doi       = {10.1016/j.spa.2017.01.001}
}

@article{LelievreZhang2019,
  author    = {T. Leli{\`e}vre and W. Zhang},
  title     = {Pathwise estimates for effective dynamics: The case of nonlinear vectorial reaction coordinates},
  journal   = {Multiscale Modeling \& Simulation},
  volume    = {17},
  year      = {2019}
}

@article{MardtPasqualiWuNoe2018,
  author    = {A. Mardt and L. Pasquali and H. Wu and F. No{\'e}},
  title     = {VAMPnets for deep learning of molecular kinetics},
  journal   = {Nature Communications},
  volume    = {9},
  number    = {5},
  year      = {2018},
  doi       = {10.1038/s41467-017-02388-1}
}

@article{WuNoe2024,
  author    = {H. Wu and F. No{\'e}},
  title     = {Reaction coordinate flows for model reduction of molecular kinetics},
  journal   = {The Journal of Chemical Physics},
  volume    = {160},
  year      = {2024},
  pages     = {044109},
  doi       = {10.1063/5.0176078}
}

@article{LuVandenEijnden2014,
  author    = {Jianfeng Lu and Eric Vanden‐Eijnden},
  title     = {Exact dynamical coarse‐graining without time‐scale separation},
  journal   = {The Journal of Chemical Physics},
  volume    = {141},
  number    = {4},
  pages     = {044109},
  year      = {2014},
  doi       = {10.1063/1.4890367}
}

@article{E2005Transition,
  title={Transition pathways in complex systems: Reaction coordinates, isocommittor surfaces, and transition tubes},
  author={E, Weinan and Ren, Weiqing and Vanden-Eijnden, Eric},
  journal={Chemical Physics Letters},
  volume={413},
  number={1-3},
  pages={242--247},
  year={2005},
  publisher={Elsevier},
  doi={10.1016/j.cplett.2005.07.084}
}

@article{PerezHernandez2013,
  author    = {Gabriel P{\'e}rez-Hern{\'a}ndez and Frank Paul and Thomas Giorgino and Gianni De Fabritiis and Frank No{\'e}},
  title     = {Identification of slow molecular order parameters for Markov model construction},
  journal   = {The Journal of Chemical Physics},
  volume    = {139},
  number    = {1},
  pages     = {015102},
  year      = {2013},
  doi       = {10.1063/1.4811489}
}

@article{DelyonHu2006,
  author  = {Delyon, Bernard and Hu, Ying},
  title   = {Simulation of conditioned diffusion and application to parameter estimation},
  journal = {Stochastic Processes and their Applications},
  volume  = {116},
  number  = {11},
  pages   = {1660--1675},
  year    = {2006},
  doi     = {10.1016/j.spa.2006.04.004}
}

@article{SchauerVanderMeulenVanZanten2017,
  author  = {Schauer, Moritz and {van der Meulen}, Frank and {van Zanten}, Harry},
  title   = {Guided proposals for simulating multi-dimensional diffusion bridges},
  journal = {Bernoulli},
  volume  = {23},
  number  = {4A},
  pages   = {2917--2950},
  year    = {2017},
  doi     = {10.3150/16-BEJ833}
}

@article{Pedersen1995,
  author  = {Pedersen, Asger Roer},
  title   = {A New Approach to Maximum Likelihood Estimation for Stochastic Differential Equations Based on Discrete Observations},
  journal = {Scandinavian Journal of Statistics},
  volume  = {22},
  number  = {1},
  pages   = {55--71},
  year    = {1995}
}

@article{Orland2011LangevinBridges,
  author  = {Orland, Henri},
  title   = {Generating transition paths by {L}angevin bridges},
  journal = {The Journal of Chemical Physics},
  year    = {2011},
  volume  = {134},
  number  = {17},
  pages   = {174114},
  doi     = {10.1063/1.3586036}
}

@misc{DelarueKoehlOrland2016ConditionedLangevin,
  author        = {Delarue, Marc and Koehl, Patrice and Orland, Henri},
  title         = {Conditioned {L}angevin {D}ynamics enables efficient sampling of transition paths},
  year          = {2016},
  eprint        = {1611.07657},
  archivePrefix = {arXiv},
  primaryClass  = {cond-mat.stat-mech},
  doi           = {10.48550/arXiv.1611.07657},
  url           = {https://arxiv.org/abs/1611.07657}
}

@article{DelarueKoehlOrland2017AbInitioCLD,
  author  = {Delarue, Marc and Koehl, Patrice and Orland, Henri},
  title   = {\textit{Ab initio} sampling of transition paths by conditioned {L}angevin dynamics},
  journal = {The Journal of Chemical Physics},
  year    = {2017},
  volume  = {147},
  number  = {15},
  pages   = {152703},
  doi     = {10.1063/1.4985651}
}

@article{DarvePohorille2001,
  author  = {Darve, Eric and Pohorille, Andrew},
  title   = {Calculating Free Energies Using Average Force},
  journal = {The Journal of Chemical Physics},
  year    = {2001},
  volume  = {115},
  number  = {20},
  pages   = {9169--9183},
  doi     = {10.1063/1.1410978}
}

@article{ComerEtAl2015ABF,
  author  = {Comer, Jeffrey and Gumbart, James C. and H{\'e}nin, J{\'e}r{\^o}me and Leli{\`e}vre, Tony and Pohorille, Andrew and Chipot, Christophe},
  title   = {The Adaptive Biasing Force Method: Everything You Always Wanted To Know but Were Afraid To Ask},
  journal = {The Journal of Physical Chemistry B},
  year    = {2015},
  volume  = {119},
  number  = {3},
  pages   = {1129--1151},
  doi     = {10.1021/jp506633n}
}

@article{TorrieValleau1977,
  author  = {Torrie, G. M. and Valleau, J. P.},
  title   = {Nonphysical Sampling Distributions in Monte Carlo Free-Energy Estimation: Umbrella Sampling},
  journal = {Journal of Computational Physics},
  year    = {1977},
  volume  = {23},
  number  = {2},
  pages   = {187--199},
  doi     = {10.1016/0021-9991(77)90121-8}
}

@article{KumarEtAl1992WHAM,
  author  = {Kumar, Shankar and Rosenberg, John M. and Bouzida, Robert and Swendsen, Robert H. and Kollman, Peter A.},
  title   = {The Weighted Histogram Analysis Method for Free-Energy Calculations on Biomolecules. I. The Method},
  journal = {Journal of Computational Chemistry},
  year    = {1992},
  volume  = {13},
  number  = {8},
  pages   = {1011--1021},
  doi     = {10.1002/jcc.540130812}
}

@article{LaioParrinello2002,
  author  = {Laio, Alessandro and Parrinello, Michele},
  title   = {Escaping Free-Energy Minima},
  journal = {Proceedings of the National Academy of Sciences},
  year    = {2002},
  volume  = {99},
  number  = {20},
  pages   = {12562--12566},
  doi     = {10.1073/pnas.202427399}
}

@article{ValssonParrinello2014VES,
  author  = {Valsson, Omar and Parrinello, Michele},
  title   = {Variational Approach to Enhanced Sampling and Free Energy Calculations},
  journal = {Physical Review Letters},
  year    = {2014},
  volume  = {113},
  number  = {9},
  pages   = {090601},
  doi     = {10.1103/PhysRevLett.113.090601}
}

@article{InvernizziParrinello2020OPES,
  author  = {Invernizzi, Michele and Parrinello, Michele},
  title   = {Rethinking Metadynamics: From Bias Potentials to Probability Distributions},
  journal = {The Journal of Physical Chemistry Letters},
  year    = {2020},
  volume  = {11},
  number  = {7},
  pages   = {2731--2736},
  doi     = {10.1021/acs.jpclett.0c00497}
}

@article{E_Ren_VandenEijnden2002String,
  author  = {E, Weinan and Ren, Weiqing and Vanden-Eijnden, Eric},
  title   = {String Method for the Study of Rare Events},
  journal = {Physical Review B},
  year    = {2002},
  volume  = {66},
  number  = {5},
  pages   = {052301},
  doi     = {10.1103/PhysRevB.66.052301}
}

@article{MaraglianoEtAl2006StringCV,
  author  = {Maragliano, Luca and Fischer, Andreas and Vanden-Eijnden, Eric and Ciccotti, Giovanni},
  title   = {String Method in Collective Variables: Minimum Free Energy Paths and Isocommittor Surfaces},
  journal = {The Journal of Chemical Physics},
  year    = {2006},
  volume  = {125},
  number  = {2},
  pages   = {024106},
  doi     = {10.1063/1.2212942}
}

@article{DellagoEtAl1998TPS,
  author  = {Dellago, Christoph and Bolhuis, Peter G. and Csajka, F{\'e}lix S. and Chandler, David},
  title   = {Transition Path Sampling and the Calculation of Rate Constants},
  journal = {The Journal of Chemical Physics},
  year    = {1998},
  volume  = {108},
  number  = {5},
  pages   = {1964--1977},
  doi     = {10.1063/1.475562}
}

@article{vanErpMoroniBolhuis2003TIS,
  author  = {van Erp, Titus S. and Moroni, Daniele and Bolhuis, Peter G.},
  title   = {A Novel Path Sampling Method for the Calculation of Rate Constants},
  journal = {The Journal of Chemical Physics},
  year    = {2003},
  volume  = {118},
  number  = {17},
  pages   = {7762--7774},
  doi     = {10.1063/1.1562614}
}

@article{AllenFrenkelTenWolde2006FFS,
  author  = {Allen, Rosalind J. and Frenkel, Daan and ten Wolde, Pieter Rein},
  title   = {Forward Flux Sampling-Type Schemes for Simulating Rare Events: Efficiency Analysis},
  journal = {The Journal of Chemical Physics},
  year    = {2006},
  volume  = {124},
  number  = {19},
  pages   = {194111},
  doi     = {10.1063/1.2198827}
}

@article{HuberKim1996WE,
  author  = {Huber, G. A. and Kim, S.},
  title   = {Weighted-Ensemble Brownian Dynamics Simulations for Protein Association Reactions},
  journal = {Biophysical Journal},
  year    = {1996},
  volume  = {70},
  number  = {1},
  pages   = {97--110},
  doi     = {10.1016/S0006-3495(96)79552-8}
}

@article{CerouGuyader2007AMS,
  author  = {C{\'e}rou, Fr{\'e}d{\'e}ric and Guyader, Arnaud},
  title   = {Adaptive Multilevel Splitting for Rare Event Analysis},
  journal = {Stochastic Analysis and Applications},
  year    = {2007},
  volume  = {25},
  number  = {2},
  pages   = {417--443},
  doi     = {10.1080/07362990601139628}
}

\appendix

\section{Technical Derivations for the Guided Bridge Lifting Method}
\label{app:guided-bridge}

\subsection{It\^{o} formula for the CV dynamics}
\label{app:ito}
Consider the controlled overdamped Langevin dynamics
\begin{equation}
dX_t = \bigl(b(X_t) + u(t,X_t)\bigr)\,dt + \sigma\, dW_t,
\label{eq:app_controlledSDE}
\end{equation}
where $W_t$ is $d$-dimensional Brownian motion and $\sigma>0$ is constant. Define
\[
Z_t := \xi(X_t)\in\mathbb{R}^m.
\]
Applying It\^{o}'s formula componentwise to each $\xi_i(X_t)$ yields
\[
d\xi_i(X_t)
= \nabla \xi_i(X_t)\cdot dX_t
+ \frac{1}{2}\sum_{p,q=1}^d (\sigma\sigma^\top)_{pq}\,\partial_{pq}\xi_i(X_t)\,dt.
\]
Since $\sigma$ is constant and isotropic, $\sigma\sigma^\top = \sigma^2 I_d$, hence
\[
d\xi_i(X_t)
= \nabla \xi_i(X_t)\cdot (b(X_t)+u(t,X_t))\,dt
+ \frac{\sigma^2}{2}\Delta \xi_i(X_t)\,dt
+ \sigma\,\nabla \xi_i(X_t)\cdot dW_t.
\]
Using $J_\xi=\nabla \xi$ and collecting components gives the vector form
\begin{equation}
dZ_t
=
J_\xi(X_t)\bigl(b(X_t)+u(t,X_t)\bigr)\,dt
+ \frac{\sigma^2}{2}\Delta \xi(X_t)\,dt
+ \sigma\,J_\xi(X_t)\,dW_t,
\label{eq:app_ito_vector}
\end{equation}
where $\Delta \xi(x) := (\Delta \xi_1(x),\ldots,\Delta \xi_m(x))^\top$.
Using the generator of the uncontrolled system, see (\ref{cL2}), this reduced to
\[
dZ_t
= (\mathcal{L}\xi)(X_t)\,dt+
J_\xi(X_t)u(t,X_t)\,dt
+ \sigma\,J_\xi(X_t)\,dW_t.
\]
This equation is unclosed in principle because the right hand side still depends on $X_t$. The equation (\ref{eq:effSDE}) for the effective dynamics results from closing it by setting $u=0$, and replacing $\mathcal{L}\xi$ as well as $J_\chi(X_t)$ by its $\mu$-weighted average along the respective fiber $\mathbb{L}_{\chi(X_t)}$.

The special form of the effective dynamics (\ref{eff-dyn}) for the specific CV $\xi=\chi$ additionally results from the special identity $\mathcal{L}\chi=(a+\lambda \chi)$.

\subsection{Computing $p(s,z)$}
The transition probability $p(s,z)$ is given via the backward Kolmogorov PDE (\ref{eq:bk-effective}) of the effective dynamics, or, respectively, via 
\[
p(s,z)=\exp((t-s)\cL_{\text{eff}})\mathbf{1}_{z>z_\ast}=\mathcal{K}^{t-s}_{\text{eff}}\mathbf{1}_{z>z_\ast}.
\]
Let us  consider a 1-dimensional effective dynamics in the subsequent. 
The eigenfunctions $(\phi_n)_{n\ge 1}$ of $\cL_{\text{eff}}$ form an eigenbasis of $L^2(\pi)$ where $\pi$ is the invariant density of the effective dynamics, given by 
\[
\pi(z)=\frac{1}{\mathcal{Z}}\exp(-V_{\text{eff}}(z)), \quad\mathcal{Z}=\int_0^1 \exp(-V_{\text{eff}}(z))dz.
\]
Let us assume that the $\phi_n$ denote $L^2_\pi$-normalized  eigenfunctions of $\cL_{\text{eff}}$, ordered such that the corresponding eigenvalues satisfy $0=\lambda_1>\lambda_2>\lambda_3>\ldots$.
For $B_{z_\ast}=(z_\ast,1]$ the spectral expansion reads
\[
p(s,z)
=
\sum_{n\ge 1} e^{\lambda_n (t-s)}\,\phi_n(z)\,
\int_{z_\ast}^1 \phi_n(\zeta)\,\pi(\zeta)\,d\zeta.
\]
With $\phi_1\equiv 1$, and assuming strong metastability with two metastable components, i.e., $\lambda_2\gg\lambda_3$, we get for all $s<t$ not too close to $t$,
\[
p(s,z)\approx \pi(B_{z_\ast}) + e^{\lambda_2(t-s)} \phi_2(z)\int_{z_\ast}^1 \phi_2(\zeta)\pi(\zeta)d\zeta. 
\]
For the one-dimensional CV $\chi$ as considered above, the drift of the effective dynamics is given by $b(z)=c+\lambda z$. In this case, $\phi_2$ and $\lambda_2$ are known:
\[
\lambda_2=\lambda,\qquad \phi_2(z)=\frac{1}{a}(c+\lambda z),\quad a^2 = \int_0^1 (c+\lambda z)^2\pi(z) dz.
\]
With the abbreviation
\[
\gamma_\ast=\int_{z_\ast}^1\phi_2(\zeta)\pi(\zeta)d\zeta=\frac{1}{a}\int_{z_\ast}^1 (c+\lambda z)\pi(z)dz,
\]
we thus have that, for $t-s$ large enough,
\begin{equation}\label{eq:SpecApprox_p}
p(s,z)\approx \pi(B_{z_\ast}) + \frac{\gamma_\ast}{a} e^{\lambda(t-s)} (c+\lambda z),
\end{equation}
and, therefore,
\begin{equation}\label{eq:derivativeApprox}
\partial_z\log p(s,z) = \frac{\partial_z p(s,z)}{p(s,z)}\approx \frac{\gamma_\ast\,\lambda e^{\lambda_2(t-s)} }{a\pi(B_{z_\ast}) + \gamma_\ast e^{\lambda(t-s)} (c+\lambda z)}.
\end{equation}

\begin{figure}
    \centering
    \includegraphics[width=0.7\linewidth]{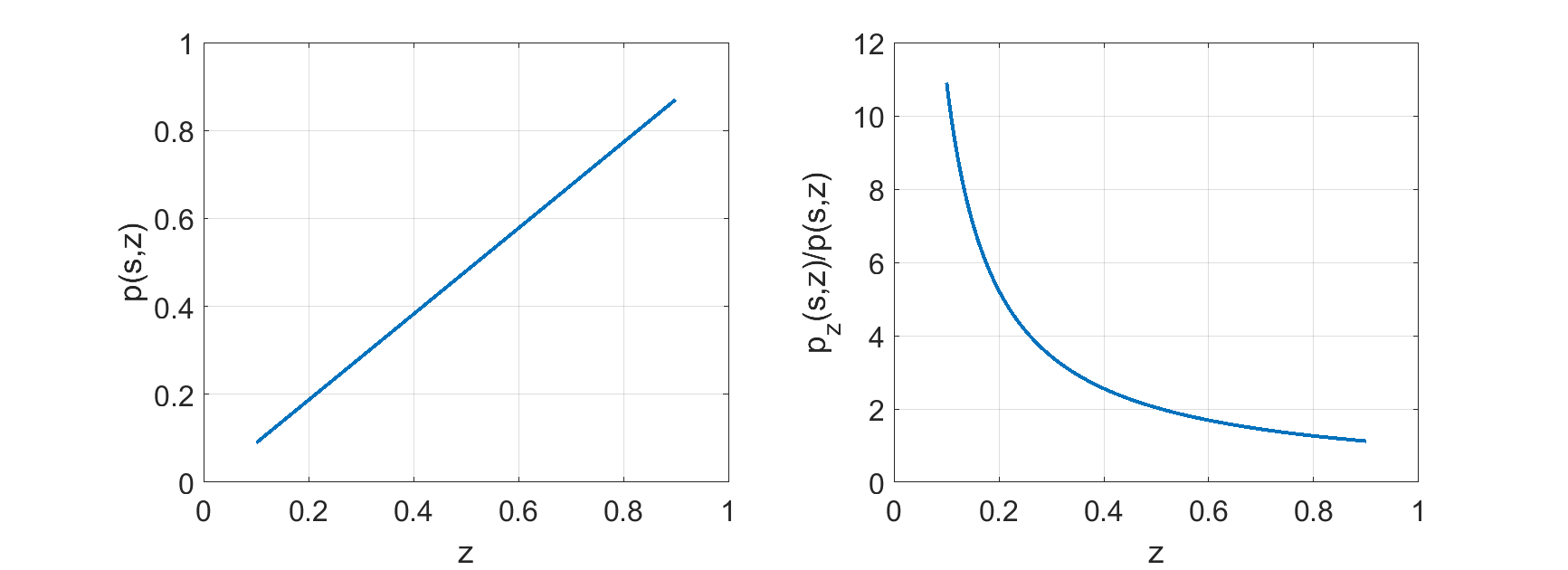}
    \caption{Spectral approximation of $p(s,z)$ (left) and of $\partial_z\log p(s,z)=p_z(s,z)/p(s,z)$ (right) for $s=0$ via equations (\ref{eq:SpecApprox_p}) and (\ref{eq:derivativeApprox}) for the simple test systems considered in Sec.~\ref{sec:numexp} and final time $t=20$.}
    \label{fig:SpectralApprox_p}
\end{figure}

\section{High-dimensional test system}\label{sec:high-d test}
We can define a high-dimensional test system based on the simple test system given in Sec.~\ref{sec:numexp}.
For the full-dimension state vector $x=(x_1,\ldots,x_d)^\top$, $d\ge 2$, we define the full-dimension potential 
\[
W(x)=V_{dw}(x_1,x_2)+\frac{1}{2}\sum_{j=3}^d \omega_j^2x_j^2.
\]
We assume that the constants $\omega_j$ are chosen such that the membership function
\[
\hat{\chi}(x)=\chi(x_1,x_2)
\]
is an appropriate CV for the system. Let $R$ denote an orthonormal $d\times d$ matrix with rows $R_j^\top$,
\[
R=\begin{pmatrix}
    R_1^\top\\ \vdots\\R_d^\top 
\end{pmatrix},\quad\text{and submatrix}\quad \mathcal{R}=\begin{pmatrix}
    R_1^\top\\R_2^\top 
\end{pmatrix}
\]
and define the full-dimensional potential via 
\[
V(x)=W(Rx).
\]
Correspondingly, we consider the CV of the system to be 
\[
\xi(x)=\hat{\chi}(Rx)=\chi(\mathcal{R} x),
\]
so that its latent space is 1-dimensional and identical to $[0,1]$. For $z\in [0,1]$, the level set of $\xi$ is
\[
\mathbb{L}_z = \{x:\;\xi(x)=z\}=\{x:\; \hat{\chi}(R x)=z\}=\{x:\;\chi(R_1^\top x,R_2^\top x)=z\}
\]

\end{document}